\documentclass[a4paper,10pt]{article}
\usepackage[T1]{fontenc}
\usepackage[utf8]{inputenc}
\usepackage{lmodern}
\usepackage[left=3.5cm,right=3.5cm,top=1.88cm,bottom=2.82cm]{geometry}
\usepackage{amsmath,amssymb,amsthm,mathtools}
\usepackage{microtype}
\usepackage{needspace}
\usepackage{enumitem}
\usepackage{tikz}
\usetikzlibrary{arrows.meta,calc,positioning,fit,backgrounds}
\usepackage[section]{placeins}
\tikzset{
	vtx/.style={circle,fill=black,inner sep=1.4pt},
	selected/.style={black,line width=.9pt},
	omitted/.style={gray!60,densely dashed,line width=.45pt},
	oldregion/.style={draw=gray!55,fill=gray!6,rounded corners=3pt}
}
\usepackage[hidelinks]{hyperref}
\hypersetup{pdftitle={Independent and end-faithful spanning trees in infinite graphs},pdfauthor={Rodrigo Santos Monteiro}}
\newtheorem{theorem}{Theorem}[section]
\AddToHook{env/theorem/before}{\Needspace{6\baselineskip}}
\newcommand{\mainresultnumber}{}
\newtheorem*{mainresult}{Theorem~\mainresultnumber}
\newenvironment{restated}[1]
{\renewcommand{\mainresultnumber}{\ref{#1}}\begin{mainresult}}
	{\end{mainresult}}
\newtheorem{claim}[theorem]{Claim}
\newtheorem{proposition}[theorem]{Proposition}
\newtheorem{conjecture}[theorem]{Conjecture}
\newtheorem{question}[theorem]{Question}
\theoremstyle{definition}
\newtheorem{definition}[theorem]{Definition}
\newtheorem{example}[theorem]{Example}
\numberwithin{equation}{section}
\numberwithin{figure}{section}

\title{Independent and end-faithful spanning trees in infinite graphs}
\author{Rodrigo Santos Monteiro}

\newcommand{\Addresses}{%
	\par\bigskip
	\begingroup
	\footnotesize
	\noindent R.~S.~Monteiro, Instituto de Ciências Matemáticas e de Computação,
	Universidade de São Paulo\\
	Avenida Trabalhador São-Carlense, 400, São Carlos, SP, 13566-590, Brazil
	\par\nopagebreak
	\noindent\textit{E-mail}: \texttt{rodrigosm@usp.br}
	\par
	\endgroup
}

\date{}

\begin{document}
	\maketitle
%	\linenumbers
	
	\begin{abstract}
		Spanning trees rooted at the same vertex are independent if, for each vertex other than the root, the paths from the root to that vertex are pairwise distinct and internally vertex-disjoint. Zehavi and Itai conjectured that, for every positive integer $k$, every finite $k$-connected graph has $k$ independent spanning trees rooted at any prescribed vertex. We show that every $2$-connected graph, finite or infinite, has two independent spanning trees rooted at any prescribed vertex. For locally finite $2$-connected graphs, both trees can be chosen end-faithful. We also prove that every locally finite $3$-connected chordal graph has three independent end-faithful spanning trees. Counterexamples show that local finiteness is essential for the version with independent end-faithful spanning trees. We further show that every $3$-connected graph with a universal vertex or a finite vertex cover has three independent spanning trees rooted at any prescribed vertex.
	\end{abstract}
	
	
	
	\section{Introduction}\label{sec:intro}
	
	All graphs in this paper are simple. We use the terminology and notation of~\cite{Diestel}. For a positive integer $k$, a graph $G$ is \emph{$k$-connected} if $|V(G)|\geq k+1$ and $G-X$ is connected for every set $X\subseteq V(G)$ with $|X|<k$. For a tree $T$ and vertices $x,y\in V(T)$, write $xTy$ for the unique $x$--$y$ path in $T$. A \emph{spanning tree} of $G$ is a connected subgraph $T$ without cycles such that $V(T)=V(G)$. Specifying a vertex $r\in V(T)$ makes $T$ a \emph{rooted tree}, with root $r$.
	
	\begin{definition}\label{def:independent}
		Let $k$ be a positive integer, let $G$ be a graph, and let $r\in V(G)$. Spanning trees $T_1,\ldots,T_k$ of $G$ are \emph{independent}, rooted at $r$, if, for every $v\in V(G)\setminus\{r\}$ and distinct $i,j\in\{1,\ldots,k\}$, the paths $rT_iv$ and $rT_jv$ are distinct and satisfy $V(rT_iv)\cap V(rT_jv)=\{r,v\}$.
	\end{definition}

	\begin{conjecture}[Zehavi and Itai~\cite{ZI}]\label{conj:finite}
		For every positive integer $k$, every finite $k$-connected graph has $k$ independent spanning trees rooted at any prescribed vertex.
	\end{conjecture}
	
	The case $k=1$ follows from the existence of a spanning tree in every finite connected graph. Itai and Rodeh~\cite[Section~7.2]{IR} proved the case $k=2$. Cheriyan and Maheshwari~\cite{CM} and, independently, Zehavi and Itai~\cite[Theorem~1.2]{ZI} proved the case $k=3$, and Curran, Lee and Yu~\cite[Theorem~1.1]{CLY} proved the case $k=4$. We record these cases together; the case $k=3$ will also be used in our proofs.
	
	\begin{theorem}[\cite{IR,CM,ZI,CLY}]\label{thm:finite}
		Let $k\in\{2,3,4\}$. Every finite $k$-connected graph has $k$ independent spanning trees rooted at any prescribed vertex.
	\end{theorem}
	
	Conjecture~\ref{conj:finite} remains open in general for $k\geq5$;
	recent progress for random graphs appears in~\cite{DFKPY,HLMPW}. The survey of Cheng, Wang and Fan~\cite{CWF} gives a broader account of independent spanning trees in finite networks and their applications.
	
	This paper studies the independent-spanning-tree problem in infinite graphs. To the best of our knowledge, the independent-spanning-tree problem has not previously been studied in infinite graphs.
	
	Our main results are Theorems~\ref{thm:two}, \ref{thm:faithful} and~\ref{thm:chordal}. The first extends the case $k=2$ for any $2$-connected graph.
	\begin{restated}{thm:two}
		Every $2$-connected graph has two independent spanning trees rooted at any prescribed vertex.
	\end{restated}
	
	Freudenthal introduced ends in 1931 in his study of topological spaces and topological groups~\cite{Freudenthal}. Independently, Halin introduced ends of graphs and end-faithful spanning trees in 1964~\cite{Halin}. A graph is \emph{locally finite} if every vertex has finite degree. For connected locally finite graphs, the topological and graph-theoretical notions of an end agree. We use the graph-theoretical notion and ask how independent spanning trees represent the ends of the graph. We recall the definitions below.
	
	\begin{definition}\label{def:ends}
		Let $G$ be a graph. A one-way infinite path is called a \emph{ray}, and its infinite connected subgraphs are its \emph{tails}. We denote the set of rays of $G$ by $\mathcal{R}(G)$. Two rays $r,s\in\mathcal{R}(G)$ are \emph{equivalent in $G$} if, for every finite set $X\subseteq V(G)$, $r$ and $s$ have tails in the same component of $G-X$. Each equivalence class $[r]_G$ is called an \emph{end} of $G$. We denote the set of ends of $G$ by $\Omega(G)$.
	\end{definition}
	
	\begin{definition}\label{def:faithful}
		A subgraph $H\subseteq G$ \emph{represents} an end $[r]_G\in\Omega(G)$ if $H$ contains a ray equivalent to $r$ in $G$. For a spanning tree $T$ of $G$, the inclusion $\mathcal{R}(T)\subseteq\mathcal{R}(G)$ induces a map $\Omega(T)\to\Omega(G)$ given by $[r]_T\mapsto[r]_G$. The tree $T$ is \emph{end-faithful} if this map is a bijection.
	\end{definition}
	
	Equivalently, a spanning tree is end-faithful if each end of $G$ contains exactly one ray of the tree starting at its root~\cite[Section~3]{DK}. 
	The existence of a single end-faithful spanning tree has been studied extensively; see Diestel~\cite{DiestelEnds}, \v{S}ir\'a\v{n}~\cite{Siran}, and the more recent work of B\"urger and Kurkofka~\cite{BK}. In particular, Seymour and Thomas~\cite{ST} constructed an infinitely connected graph with no end-faithful spanning tree, disproving the assertion that every connected graph admits one. 
	
	For locally finite graphs, we strengthen Theorem~\ref{thm:two} as follows.
	
	\begin{restated}{thm:faithful}
		Every locally finite $2$-connected graph has two independent end-faithful spanning trees rooted at any prescribed vertex.
	\end{restated}
	
	\begin{definition}\label{def:chordal}
		A \emph{clique} of a graph is a set of pairwise adjacent vertices; it is \emph{maximal} if it is not properly contained in another clique. A graph is \emph{chordal} if it has no induced cycle of length at least four, where a cycle is induced when no two of its nonconsecutive vertices are adjacent.
	\end{definition}
	
	Chordal graphs are an important class in structural and algorithmic graph theory: finite chordal graphs admit efficient algorithms for minimum colouring, maximum clique and maximum independent set~\cite{Gavril72}. Their structure can be described by clique trees, which organise their maximal cliques and the separators formed by intersections of adjacent cliques~\cite{Gavril74}. Hofer-Temmel and Lehner extended this representation to connected locally finite chordal graphs~\cite{HTL}, as stated in Theorem~\ref{thm:cliquetree}. For this class, we obtain the following result.
	
	\begin{restated}{thm:chordal}
		Every locally finite $3$-connected chordal graph has three independent end-faithful spanning trees rooted at any prescribed vertex.
	\end{restated}
	
	Examples~\ref{ex:fan} and~\ref{ex:squarefan} show that local finiteness is essential in Theorems~\ref{thm:faithful} and~\ref{thm:chordal}, respectively.
	For $k=3$, we also prove a partial extension of Theorem~\ref{thm:finite} to infinite graphs. A vertex is \emph{universal} if it is adjacent to every other vertex, and a \emph{vertex cover} is a set of vertices meeting every edge.
	
	\begin{restated}{thm:three}
		Every $3$-connected graph with a universal vertex or a finite vertex cover has three independent spanning trees rooted at any prescribed vertex.
	\end{restated}
	
	Question~\ref{ques:three} asks whether this conclusion holds for every infinite $3$-connected graph. In the locally finite case, we also ask whether the three trees can be chosen end-faithful (Question~\ref{ques:threefaithful}).

	\section{Independent spanning trees}\label{sec:ordered}
	
	For a set $X$ of vertices, the graph $G-X$ is obtained by deleting $X\cap V(G)$ and all edges incident with these vertices. We write $G-x$ for $G-\{x\}$ and $G-H$ for $G-V(H)$ when $H$ is a subgraph of $G$. For $W\subseteq V(G)$, the graph $G[W]$ has vertex set $W$ and all edges of $G$ with both ends in $W$. For $F\subseteq E(G)$, the graph $G-F$ is obtained by deleting the edges in $F$ and retaining all vertices. We also write $G-e$ for deletion of an edge and $G+e$ for addition of an edge, including its endvertices if necessary. For $U\subseteq V(G)$, the set $N_G(U)$ consists of the vertices outside $U$ with a neighbour in $U$; we put $N_G(v)=N_G(\{v\})$ and $N_G(H)=N_G(V(H))$ for a subgraph $H$. The distance $d_G(x,y)$ between two vertices in the same component is the minimum number of edges of an $x$--$y$ path. The ball of radius $n$ around $x$ is the set of vertices $y$ in the component of $x$ with $d_G(x,y)\leq n$. For a subgraph $H$ of $G$, an \emph{$H$-path} is a path in $G$ with distinct endvertices in $H$, all internal vertices outside $H$, and no edge in $H$. Every $H$-path used to extend a subgraph below has length at least two.
	
	\begin{definition}\label{def:rooted}
		Let $T$ be a tree rooted at $r$. The \emph{depth} of $v\in V(T)$ is the number of edges of $rTv$. A vertex $x$ is an \emph{ancestor} of $y$ if $x\in V(rTy)$; then $y$ is a \emph{descendant} of $x$. The ancestor is \emph{proper} if $x\neq y$. For $v\neq r$, the \emph{parent} of $v$ is its predecessor on $rTv$, and the \emph{children} of a vertex are the vertices having it as their parent. The subgraph induced in $T$ by a vertex and all its descendants is its \emph{descendant subtree}.
	\end{definition}
	
	A \emph{leaf} is a vertex of degree one in a tree. By attaching a new vertex $v$ as a leaf to $u$, we mean adding $v$ and the edge $uv$ to the tree.
	We use the following form of Menger's theorem.
	
	\begin{theorem}[Menger's theorem; see~{\cite[Corollary~3.3.4 and Section~8.4]{Diestel}}]\label{thm:fan}
		Let $G$ be a graph, let $k\geq1$ be an integer, let $x\in V(G)$, and let $S\subseteq V(G)\setminus\{x\}$ with $|S|\geq k$. Suppose that, for every set $X\subseteq V(G)\setminus\{x\}$ with $|X|<k$, the component of $G-X$ containing $x$ meets $S\setminus X$. Then there are $k$ paths from $x$ to distinct vertices of $S$ such that each path meets $S$ only at its last vertex and any two of the paths meet only at $x$.
	\end{theorem}
	
	In particular, if $G$ is $2$-connected, $H\subseteq G$ has at least two vertices, and $x\notin V(H)$, then $x$ lies on an $H$-path. Indeed, for $|X|<2$ with $x\notin X$, the graph $G-X$ is connected and $V(H)\setminus X$ is nonempty. Apply Theorem~\ref{thm:fan} with $S=V(H)$ and $k=2$, and concatenate its two paths at $x$.
	
	The construction for any vertex set also uses Zorn's lemma.
	
	\begin{theorem}[Zorn's lemma; see~{\cite[Appendix~A]{Diestel}}]\label{thm:zorn}
		Every nonempty partially ordered set in which every chain has an upper bound has a maximal element.
	\end{theorem}
	
	\begin{theorem}\label{thm:two}
		Every $2$-connected graph has two independent spanning trees rooted at any prescribed vertex.
	\end{theorem}
	
	\begin{proof}
		If $G$ is finite, apply Theorem~\ref{thm:finite} with $k=2$. Suppose that $G$ is infinite. We extend the construction for finite graphs of~\cite[Section~3]{ZI} to any vertex set, retaining two trees and a linear order throughout. Here $<$ denotes a strict linear order: it is irreflexive and transitive, and, for distinct $x,y$, exactly one of $x<y$ and $y<x$ holds. A path $x_0\ldots x_m$, read in the indicated direction, is \emph{increasing} if $x_0<\cdots<x_m$. A path of length zero is increasing.
		
		Consider a subgraph $H\subseteq G$, distinct vertices $s,t\in V(H)$ with $st\in E(H)$, trees $A,B$ spanning $H$, and a linear order $<$ on $V(H)$ with minimum $s$ and maximum $t$. We require:
		\begin{enumerate}[label=(\roman*),ref=(\roman*)]
			\item\label{cond:increasingA} The path $sAv$ is increasing for every $v\in V(H)$.
			\item\label{cond:increasingB} The vertex $s$ is a leaf of $B$ incident with $st$, and the path $vBt$ is increasing for every $v\in V(H)\setminus\{s\}$.
		\end{enumerate}
		Condition~\ref{cond:increasingA} implies that $t$ is a leaf of $A$. Otherwise, besides its predecessor on $sAt$, it would have a neighbour $w$ such that $sAw$ passes through $t$, contrary to the maximality of $t$ in the order.
		
		Suppose also that $st\notin E(A)$. For $v\notin\{s,t\}$, every internal vertex of $sAv$ is smaller than $v$. By condition~\ref{cond:increasingB}, the path $sBv$ consists of the edge $st$ followed by $tBv$, and all its internal vertices are larger than $v$. Thus the two root paths have disjoint internal vertex sets. They are distinct because $sBv$ contains $t$ whereas $sAv$ does not. For $v=t$, the path in $B$ is the edge $st$, whereas the path in $A$ avoids that edge. Consequently $A,B$ are independent rooted at $s$ whenever $st\notin E(A)$.
		
		\begin{claim}\label{cl:extension}
			The trees $A,B$ and the order $<$ have extensions to trees spanning $G$ and a linear order on $V(G)$ with the same minimum $s$ and maximum $t$, satisfying conditions~\ref{cond:increasingA} and~\ref{cond:increasingB}. The extensions preserve both trees as subgraphs and preserve the order on $V(H)$.
		\end{claim}
		
		\begin{proof}
			We first give an extension that adds the internal vertices of one path. Suppose $A',B'$ and $<'$ are extensions satisfying conditions~\ref{cond:increasingA} and~\ref{cond:increasingB} on a vertex set $W$, and let $P=a x_1\ldots x_q b$ be a $G[W]$-path with $q\geq1$. Reverse its notation if necessary so that $a<'b$. Extend $<'$ by putting $x_1,\ldots,x_q$, in this order, immediately after $a$. Every old vertex at most $a$ precedes them, and they precede every old vertex larger than $a$, in particular $b$. Extend the trees to
			\begin{equation}\label{eq:ear}
				A'\cup a x_1\ldots x_q,
				\qquad
				B'\cup b x_q\ldots x_1.
			\end{equation}
			Each added path meets its old tree only at its attachment vertex, namely $a$ in the first tree and $b$ in the second. Hence both resulting graphs are connected and acyclic. Since they contain the old trees, uniqueness of paths in a tree implies that every path between old vertices remains unchanged.
			
			For $1\leq j\leq q$, the new path from $s$ to $x_j$ in the first tree is $sA'a$ followed by $a x_1\ldots x_j$, and the new path from $x_j$ to $t$ in the second is $x_j\ldots x_q b$ followed by $bB't$. Both are increasing. Since the new vertices lie between $a$ and $b$, they lie strictly between $s$ and $t$. Moreover, $b\neq s$, so no edge is added at $s$ in the second tree. Thus the extrema of the order and conditions~\ref{cond:increasingA} and~\ref{cond:increasingB} are preserved.
			
			Figure~\ref{fig:ear} shows the two extensions in~\eqref{eq:ear}.
			
			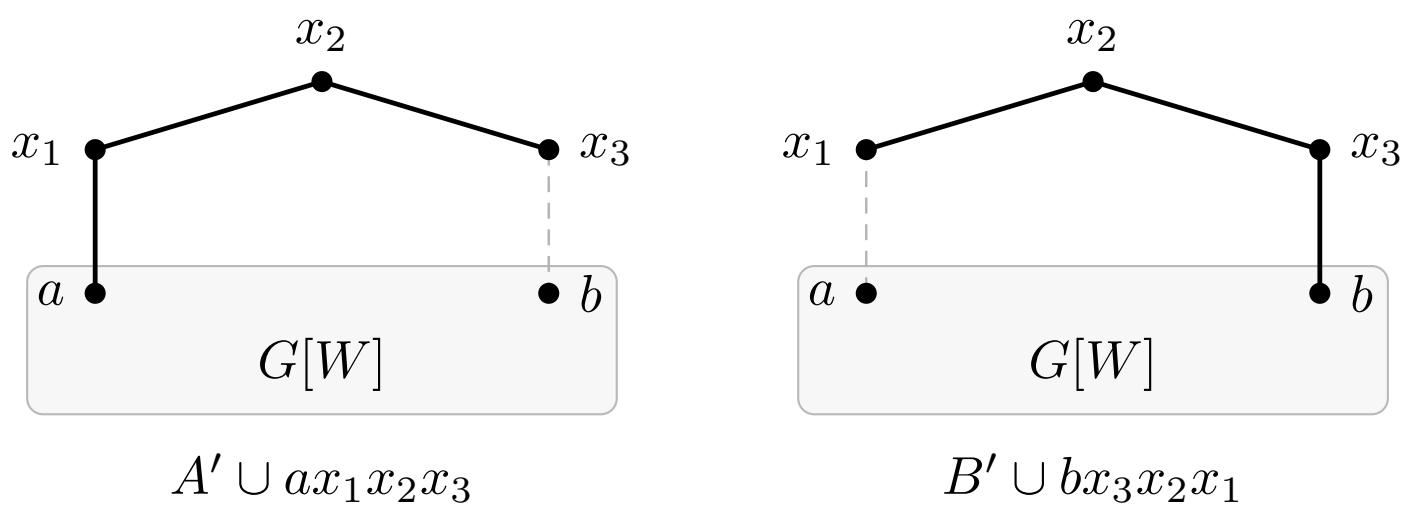
\begin{figure}[htbp]
				\centering
				\begin{tikzpicture}[x=1cm,y=1cm,every node/.style={font=\normalsize}]
					\foreach \shift/\which in {0/1,5.1/2}{
						\begin{scope}[xshift=\shift cm]
							\draw[oldregion] (-.45,-.8) rectangle (3.45,.18);
							\node at (1.5,-.47) {$G[W]$};
							\coordinate (a) at (0,0);
							\coordinate (b) at (3,0);
							\coordinate (x1) at (0,.95);
							\coordinate (x2) at (1.5,1.4);
							\coordinate (x3) at (3,.95);
							\draw[selected] (x1)--(x2)--(x3);
							\ifnum\which=1
							\draw[selected] (a)--(x1);
							\draw[omitted] (x3)--(b);
							\node at (1.5,-1.22) {$A'\cup ax_1x_2x_3$};
							\else
							\draw[omitted] (a)--(x1);
							\draw[selected] (x3)--(b);
							\node at (1.5,-1.22) {$B'\cup bx_3x_2x_1$};
							\fi
							\node[vtx,label=left:$a$] at (a) {};
							\node[vtx,label=right:$b$] at (b) {};
							\node[vtx,label=left:$x_1$] at (x1) {};
							\node[vtx,label=above:$x_2$] at (x2) {};
							\node[vtx,label=right:$x_3$] at (x3) {};
						\end{scope}
					}
				\end{tikzpicture}
				\caption{The extension in~\eqref{eq:ear} along $ax_1x_2x_3b$. Black edges are added to the indicated tree; the dashed edge is omitted. The trees on $W$ are preserved and not shown.}
				\label{fig:ear}
			\end{figure}
			
			Let $\mathcal P$ be the set of all quadruples $(W,\prec,A_W,B_W)$ such that $V(H)\subseteq W\subseteq V(G)$, the graphs $A_W,B_W$ are trees in $G[W]$ spanning $W$ and containing $A,B$, respectively, and $\prec$ is a linear order on $W$ extending $<$. Require that its minimum and maximum are $s,t$, and that the two trees and this order satisfy conditions~\ref{cond:increasingA} and~\ref{cond:increasingB}. Write $(W,\prec,A_W,B_W)\leq (W',\prec',A_{W'},B_{W'})$ when $W\subseteq W'$, $\prec\,\subseteq\,\prec'$, $A_W\subseteq A_{W'}$, and $B_W\subseteq B_{W'}$. Here each order is viewed as a set of ordered pairs, and graph inclusion requires inclusion of vertex and edge sets. This is a partial order. The quadruple $(V(H),<,A,B)$ belongs to $\mathcal P$, so $\mathcal P$ is nonempty.
			
			Let $\mathcal C$ be a nonempty chain in $\mathcal P$. Take the union of its vertex sets, the union of its order relations, and the union of its first and second trees, respectively. Write these unions as $W_{\mathcal C},\prec_{\mathcal C},A_{\mathcal C},B_{\mathcal C}$. Every nonempty finite subset of $\mathcal C$ has a greatest element. Thus any two vertices of $W_{\mathcal C}$ occur together in one member and are joined there by a path in each tree. The two union graphs are therefore connected. A cycle in either union graph has finitely many edges, all of which would lie in a single member of $\mathcal C$, contradicting that its corresponding graph is a tree. Hence both unions are trees.
			
			For distinct vertices $x,y\in W_{\mathcal C}$, a member containing both specifies exactly one of $x\prec_{\mathcal C}y$ and $y\prec_{\mathcal C}x$. No later member can reverse that comparison. If $x\prec_{\mathcal C}y$ and $y\prec_{\mathcal C}z$, a member containing these two comparisons also gives $x\prec_{\mathcal C}z$. Thus $\prec_{\mathcal C}$ is a linear order extending $<$, with minimum $s$ and maximum $t$.
			
			For a vertex $v\in W_{\mathcal C}$, choose a member containing $v$. Its paths from $s$ in the first tree and to $t$ in the second lie in the union trees; uniqueness of tree paths shows that they are the corresponding paths in the unions. Their comparisons are preserved, so they remain increasing. Every edge at $s$ in $B_{\mathcal C}$ belongs to a member of the chain, where the only such edge is $st$. Hence condition~\ref{cond:increasingB}, including the leaf requirement, holds in the union, as does condition~\ref{cond:increasingA}. The union quadruple belongs to $\mathcal P$ and is an upper bound for $\mathcal C$. The empty chain has the initial quadruple as an upper bound.
			
			By Theorem~\ref{thm:zorn}, the partially ordered set $\mathcal P$ has a maximal element $(W,\prec,A_W,B_W)$. Suppose $v\in V(G)\setminus W$. For every $X\subseteq V(G)\setminus\{v\}$ with $|X|<2$, the graph $G-X$ is connected and $W\setminus X$ is nonempty, because $s,t\in W$. Theorem~\ref{thm:fan}, with $k=2$ and target set $W$, gives two paths from $v$ to distinct vertices of $W$, meeting only at $v$ and meeting $W$ only at their endvertices. Their union is a $G[W]$-path through $v$. The extension in~\eqref{eq:ear} gives a strictly larger member of $\mathcal P$, a contradiction. Therefore $W=V(G)$, proving Claim~\ref{cl:extension}.
		\end{proof}
		
		Now let $s$ be the prescribed root. Since $G$ is $2$-connected, $s$ has two distinct neighbours, say $v_1$ and $t$. Indeed, $s$ cannot be isolated, and if it had only one neighbour, deleting that vertex would disconnect $G$. The connected graph $G-s$ contains a finite $v_1$--$t$ path. Together with $sv_1$ and $st$, this gives a cycle $C=s v_1\ldots v_m t s$.
		Set $A=C-st$, $B=C-sv_1$ and order the vertices by $s<v_1<\cdots<v_m<t$. These trees and this order satisfy conditions~\ref{cond:increasingA} and~\ref{cond:increasingB}. Claim~\ref{cl:extension} gives spanning extensions $T_1,T_2$. The first cannot contain $st$: together with the path $C-st\subseteq T_1$, that edge would form a cycle. The separation of the root paths by the order, proved above, now shows that $T_1,T_2$ are independent.
		
		The same argument permits any prescribed initial cycle of this form. Every initial tree path is preserved because the final graph is a tree containing it. In particular, $t$ is a leaf of $T_1$ with its original incident edge, and $s$ is a leaf of $T_2$ incident with $st$. Also, the argument for unions of chains applies to a nested sequence of finite extensions without a further application of Theorem~\ref{thm:zorn}. We shall use these facts in the subsequent constructions.
	\end{proof}
	
	\begin{theorem}\label{thm:three}
		Every $3$-connected graph with a universal vertex or a finite vertex cover has three independent spanning trees rooted at any prescribed vertex.
	\end{theorem}
	
	\begin{proof}
		If $G$ is finite, apply Theorem~\ref{thm:finite} with $k=3$. Suppose that $G$ is infinite, and fix the prescribed root $r$. First suppose that $G$ has a finite vertex cover $C$. Its complement $I=V(G)\setminus C$ is an independent set, so every neighbour of a vertex of $I$ lies in $C$. For each $S\subseteq C$, let $I_S=\{x\in I:N_G(x)=S\}$. The nonempty sets $I_S$ partition $I$. Choose a set $J_S\subseteq I_S$ of three vertices if $|I_S|\geq3$, and put $J_S=I_S$ otherwise. If $r\in I_S$, make the choice so that $r\in J_S$. Let $H$ be the subgraph of $G$ induced by $C\cup\bigcup_{S\subseteq C}J_S$.
		The graph $H$ is finite, since $C$ has only finitely many subsets and each $J_S$ has at most three vertices. It contains $r$.
		
		We verify that $H$ is $3$-connected. Let $X\subseteq V(H)$ with $|X|<3$, and let $a,b\in V(H)\setminus X$. There is a finite $a$--$b$ path in $G-X$. Every vertex $x$ of this path outside $H$ belongs to some class $I_S$ with $|J_S|=3$. Choose $w_x\in J_S\setminus X$, which is possible because $|X|<3$. Replace each such occurrence of $x$ by $w_x$. Its predecessor and successor on the original path lie in $C$, and both are adjacent to $w_x$, since $N_G(w_x)=N_G(x)$. No two omitted vertices are consecutive because $I$ is independent. The replacement therefore produces a finite $a$--$b$ walk in $H-X$, possibly with repeated vertices. Removing closed subwalks gives an $a$--$b$ path. Thus $H-X$ is connected.
		
		It remains to check that $|V(H)|\geq4$. If $H=G$, this follows from the definition of $3$-connectivity. Otherwise choose $x\in V(G)\setminus V(H)$. Every vertex of a $3$-connected graph has at least three neighbours: deleting a smaller neighbourhood would leave that vertex isolated and at least one other vertex present. Hence $x$ has at least three neighbours in $C$, and $H$ contains these neighbours together with a vertex of $J_{N_G(x)}$ outside $C$. This proves the required size bound and hence $3$-connectivity of $H$.
		
		Apply Theorem~\ref{thm:finite} with $k=3$ to obtain independent spanning trees $A_1,A_2,A_3$ of $H$ rooted at $r$. For every class with $I_S\setminus J_S\neq\varnothing$, choose $w_S\in J_S\setminus\{r\}$. There are three vertices in $J_S$, so this choice is possible. Let $p_i(w_S)$ be the parent of $w_S$ in $A_i$. It belongs to $S$, because every neighbour of $w_S$ lies in $C$ and $N_G(w_S)=S$. For $i\in\{1,2,3\}$, form $T_i$ by attaching every omitted vertex $x\in I_S\setminus J_S$ as a leaf to $p_i(w_S)$.
		
		The resulting graph $T_i$ is connected: every new vertex is adjacent to a vertex of the connected tree $A_i$. It contains no cycle, since its new vertices all have degree one and $A_i$ is acyclic. It is therefore a spanning tree of $G$. The paths from $r$ to vertices of $H$ are unchanged. For an omitted vertex $x\in I_S$, the path $rT_ix$ is obtained from $rA_iw_S$ by replacing its final vertex $w_S$ by $x$. These replacements preserve the internal vertex sets, and replacing $x$ by $w_S$ recovers each original path. Thus the three new paths are distinct and pairwise internally disjoint. It follows that $T_1,T_2,T_3$ are independent.
		
		Figure~\ref{fig:finitecover} illustrates the attachment for an index $i\in\{1,2,3\}$.
		
		\begin{figure}[htbp]
			\centering
			\begin{tikzpicture}[
				x=1cm,y=1cm,
				every node/.style={font=\small},
				vertex/.style={circle,fill=black,inner sep=1.4pt}
				]
				% The finite subgraph H
				\filldraw[
				fill=gray!6,
				draw=gray!55,
				rounded corners=3pt
				] (-.65,-1.95) rectangle (3.65,1.15);
				\node at (-.30,.85) {$H$};
				
				\coordinate (r) at (0,0);
				\coordinate (p) at (2.70,0);
				\coordinate (w) at (2.70,-1.25);
				\coordinate (x) at (5.30,0);
				
				% Existing path and edge
				\draw[line width=.5pt]
				(r) .. controls (.75,.22) and (1.95,.22) .. (p);
				\node at (1.35,.65) {$rA_ip_i$};
				\draw[line width=.5pt] (p)--(w);
				
				% Added edge
				\draw[line width=1.1pt] (p)--(x);
				
				% Vertices and labels
				\node[vertex,label=below:{$r$}] at (r) {};
				\node[vertex,label=above:{$p_i$}] at (p) {};
				\node[vertex,label=below:{$w_S$}] at (w) {};
				\node[vertex,label=right:{$x$}] at (x) {};
			\end{tikzpicture}
			\caption{Attaching $x$ to $p_i=p_i(w_S)$, for
				$i\in\{1,2,3\}$. The drawing shows $p_i\neq r$. Only the thick edge is added. If $p_i=r$, this edge is $rx$.}
			\label{fig:finitecover}
		\end{figure}
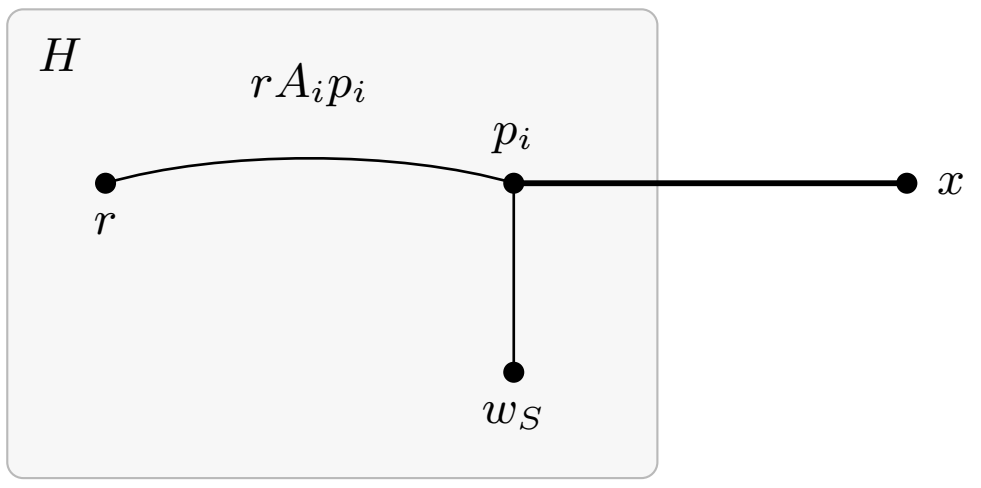
		
		Now suppose that $G$ has a universal vertex $z$, and put $F=G-z$. We have $|V(F)|\geq3$, and $F-X=G-(X\cup\{z\})$ is connected whenever $X\subseteq V(F)$ and $|X|<2$. Consequently $F$ is $2$-connected.
		
		Suppose first that $r\in V(F)$. Choose a cycle $r u\ldots t r$ in $F$, with $u\neq t$, as in the last part of the proof of Theorem~\ref{thm:two}. That proof gives independent spanning trees $A,B$ of $F$ rooted at $r$, preserving the initial paths, with $ru\in E(A)$ and with $r$ a leaf of $B$ incident with $rt$. Form $T_1$ by attaching $z$ to $u$ in $A$, form $T_2$ by attaching $z$ to $t$ in $B$, and let $T_3$ be the spanning star with centre $z$, whose edge set is $\{zv:v\in V(F)\}$. The new edges exist because $z$ is universal, and all three graphs are spanning trees of $G$.
		
		For $v\in V(F)\setminus\{r\}$, the paths $rT_1v$ and $rT_2v$ are the old independent paths in $F$, while $rT_3v=rzv$ has only $z$ as an internal vertex. They are pairwise internally disjoint and distinct. For $v=z$, the three paths are $ruz$, $rtz$ and $rz$. Since $u\neq t$, these too are distinct and pairwise internally disjoint.
		
		Finally suppose $r=z$. Choose a cycle $C=s u\ldots w t s$ in $F$, allowing $u=w$ when the cycle is a triangle. In particular, $s,t,u$ are distinct and $w\notin\{s,t\}$. Apply Claim~\ref{cl:extension} in $F$ to the trees $C-st$, $C-su$ and the order along $C-st$ from $s$ to $t$. It gives trees $A,B$ spanning $F$ and an order satisfying conditions~\ref{cond:increasingA} and~\ref{cond:increasingB}. The vertex $t$ is a leaf of $A$ incident with $wt$, and $s$ is a leaf of $B$ incident with $st$. Define spanning subgraphs of $G$ by
		\begin{equation}\label{eq:three}
			\begin{aligned}
				T_1&=(A-wt)+st+rs,\\
				T_2&=B+rt,\\
				E(T_3)&=\{rv:v\in V(F)\setminus\{s,t\}\}\cup\{su,tw\}.
			\end{aligned}
		\end{equation}
		Deleting $wt$ isolates the leaf $t$ in $A$, and adding $st$ attaches it to $s$; the further edge $rs$ attaches the new vertex $r$. Hence $T_1$ is a tree. The graph $T_2$ is a tree obtained by attaching $r$ to $t$ in $B$. The graph $T_3$ is a star on $\{r\}\cup(V(F)\setminus\{s,t\})$, with $s$ attached to $u$ and $t$ attached to $w$. This remains a tree when $u=w$. All three span $G$.
		
		Figure~\ref{fig:universalroot} shows the construction in~\eqref{eq:three} on an initial cycle of length four.
		
		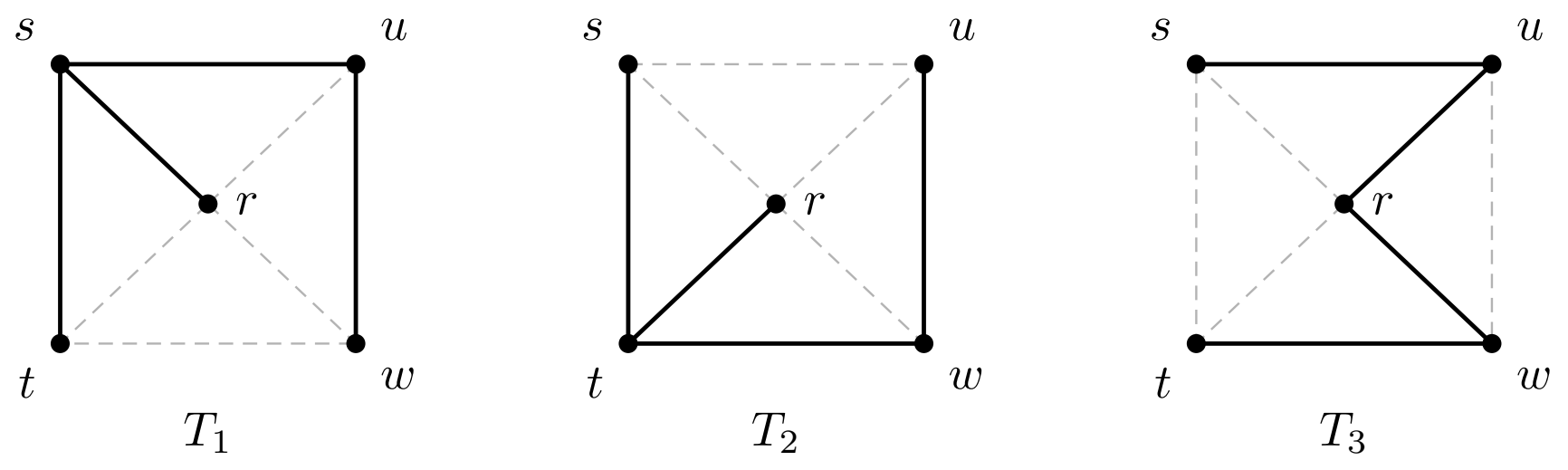
\begin{figure}[htbp]
			\centering
			\begin{tikzpicture}[x=1cm,y=1cm,every node/.style={font=\normalsize}]
				\foreach \shift/\which in {0/1,4.15/2,8.3/3}{
					\begin{scope}[xshift=\shift cm]
						\coordinate (s) at (-1.08,1.02);
						\coordinate (u) at (1.08,1.02);
						\coordinate (w) at (1.08,-1.02);
						\coordinate (t) at (-1.08,-1.02);
						\coordinate (r) at (0,0);
						\draw[omitted] (s)--(u)--(w)--(t)--(s);
						\foreach \v in {s,u,w,t}{\draw[omitted] (r)--(\v);}
						\ifnum\which=1
						\draw[selected] (r)--(s)--(u)--(w);
						\draw[selected] (s)--(t);
						\fi
						\ifnum\which=2
						\draw[selected] (s)--(t)--(w)--(u);
						\draw[selected] (r)--(t);
						\fi
						\ifnum\which=3
						\draw[selected] (s)--(u)--(r)--(w)--(t);
						\fi
						\node[vtx,label=above left:$s$] at (s) {};
						\node[vtx,label=above right:$u$] at (u) {};
						\node[vtx,label=below right:$w$] at (w) {};
						\node[vtx,label=below left:$t$] at (t) {};
						\node[vtx,label=right:$r$] at (r) {};
						\node at (0,-1.67) {$T_{\which}$};
					\end{scope}
				}
			\end{tikzpicture}
			\caption{The construction in~\eqref{eq:three} on the initial cycle $C=suwts$, with $A\cap C=C-st$ and $B\cap C=C-su$. In each panel the black edges belong to $T_i$ and the dashed edges are omitted from $T_i$. The common root $r$ is universal.}
			\label{fig:universalroot}
		\end{figure}
		
		Let $v\in V(F)\setminus\{s,t\}$. Moving the leaf $t$ in the first tree does not change $sAv$, so its root path is the edge $rs$ followed by $sAv$. Its internal vertices all belong to $F$ and are smaller than $v$, including $s$. The second root path is the edge $rt$ followed by $tBv$, and all its internal vertices in $F$ are larger than $v$, including $t$. These paths have disjoint internal vertex sets. They are distinct from one another because their first vertices after $r$ are different, and both have length at least two. The third path is the edge $rv$, so all three are distinct and pairwise internally disjoint. For the target $s$, the three paths are $rs$, $rts$ and $rus$, whose internal vertex sets are $\varnothing$, $\{t\}$ and $\{u\}$. For the target $t$, they are $rst$, $rt$ and $rwt$, whose internal vertex sets are $\{s\}$, $\varnothing$ and $\{w\}$. The stated distinctness of the cycle vertices verifies independence at these final two targets as well.
	\end{proof}
	
	\section{Two end-faithful spanning trees}\label{sec:ends}
	
	We use the standard facts that inclusion induces a function from the ends of a spanning tree to the ends of the graph and that every end of a rooted tree contains exactly one ray starting at its root~\cite[Section~3]{DK}. Thus the function in Definition~\ref{def:faithful} is a bijection  when every end of $G$ contains exactly one ray of $T$ starting at the root.
	
	For locally finite graphs, representation of every end is automatic.
	
	\begin{proposition}[Diestel and K\"uhn; see~{\cite[proof of Lemma~3.2]{DK}}]\label{prop:endrepresentation}
		Let $G$ be a connected locally finite graph, let $T$ be a spanning tree of $G$, and let $r\in V(G)$. Every end of $G$ contains a ray of $T$ starting at $r$.
	\end{proposition}
	
	We shall also use the block structure of a connected graph. A vertex $x$ of a connected graph $F$ is a \emph{cutvertex} if $F-x$ is disconnected. An edge $e$ is a \emph{bridge} if $F-e$ is disconnected. A \emph{block} of $F$ is a maximal connected subgraph with no cutvertex. For a connected graph with at least two vertices, the \emph{block-cutvertex tree} is the bipartite graph whose vertex classes are the blocks and the cutvertices of $F$, with a block $K$ adjacent to a cutvertex $x$  when $x\in V(K)$.
	
	\begin{theorem}[See~{\cite[Section~3.1, especially Lemma~3.1.4]{Diestel}}]\label{thm:blocks}
		Let $F$ be a connected graph with at least two vertices. Every edge of $F$ belongs to exactly one block. Every block is either a bridge together with its endvertices or a $2$-connected subgraph. Two distinct blocks have at most one vertex in common, and any common vertex is a cutvertex of $F$. The block-cutvertex incidence graph is a tree.
		
		For a cutvertex $x$, its deletion from this incidence tree corresponds to deletion of $x$ from $F$. For each component $\mathcal U$ of the incidence tree after deleting $x$, take the union of $V(K)\setminus\{x\}$ over the blocks $K$ in $\mathcal U$. These unions are the vertex sets of the components of $F-x$.
	\end{theorem}
	
	An \emph{edge cut} of $G$ associated with a partition $V(G)=U\mathbin{\dot\cup}W$ into two nonempty sets is the set of edges with one end in $U$ and the other in $W$. For a spanning tree $T$ of $G$ and an edge $e\in E(T)$, the \emph{fundamental cut} of $e$ is the set of all edges of $G$ with ends in the two different components of $T-e$.
	
	We also use the fundamental-cut formulation of the end-faithfulness criterion for locally finite graphs.
	
	\begin{proposition}[Diestel and K\"uhn; see~{\cite[Lemma~3.3 and the remark immediately following it]{DK}}]\label{prop:finitecutcriterion}
		Let $G$ be a connected locally finite graph and let $T$ be a spanning tree of $G$. Then $T$ is end-faithful if and only if the fundamental cut of every edge of $T$ is finite.
	\end{proposition}
	
	\begin{theorem}\label{thm:faithful}
		Every locally finite $2$-connected graph has two independent end-faithful spanning trees rooted at any prescribed vertex.
	\end{theorem}
	
	\begin{proof}
		If $G$ is finite, apply Theorem~\ref{thm:finite} with $k=2$; end-faithfulness is vacuous. Suppose that $G$ is infinite. We use the trees and linear order constructed in the proof of Theorem~\ref{thm:two}. In addition to conditions~\ref{cond:increasingA} and~\ref{cond:increasingB}, we shall ensure that every edge in either final tree has a finite fundamental cut. Proposition~\ref{prop:finitecutcriterion} will then give end-faithfulness.
		
		\begin{claim}\label{cl:finiteblocks}
			Let $H\subseteq G$ be a finite connected subgraph with at least two vertices. Every infinite component of $G-V(H)$ has finitely many blocks and finitely many cutvertices.
		\end{claim}
		
		\begin{proof}
			Let $D$ be an infinite component of $G-V(H)$, and put $Z=N_G(V(H))\cap V(D)$.
			Since $G$ is connected, $D$ has an edge to $H$, so $Z\neq\varnothing$. Local finiteness and the finiteness of $H$ imply that $Z$ is finite.
			
			Let $x$ be a cutvertex of $D$. Every component $U$ of $D-x$ meets $Z$. Otherwise $U$ would have no neighbour in $H$. It has no neighbour in another component of $G-V(H)$, and its only possible neighbour in $D-U$ is $x$. Hence $U$ would be a component of $G-x$ disjoint from the nonempty subgraph $H$, contradicting the $2$-connectivity of $G$.
			
			Choose a finite connected subgraph $L$ of $D$ containing $Z$: fix one vertex of $Z$ and join it to every other vertex of $Z$ by a path in $D$. Their union is finite and connected. If $x$ is a cutvertex of $D$, choose vertices of $Z$ in two different components of $D-x$, which is possible by the preceding paragraph. A path between them in $L$ must pass through $x$. Consequently every cutvertex of $D$ belongs to the finite vertex set $V(L)$.
			
			Write $\mathcal B$ for the block-cutvertex tree of $D$ from Theorem~\ref{thm:blocks}. Every cutvertex $x$ belongs to only finitely many blocks. Indeed, each such block contains an edge incident with $x$, distinct blocks contain distinct edges by Theorem~\ref{thm:blocks}, and $x$ has finite degree. If $D$ has no cutvertex, $\mathcal B$ has just one vertex, corresponding to the block $D$. Otherwise every block of $D$ contains a cutvertex, because $\mathcal B$ is connected. There are finitely many cutvertices, and each belongs to finitely many blocks, so $D$ has finitely many blocks. Thus $\mathcal B$ is finite.
		\end{proof}
		
		\begin{claim}\label{cl:finitecompletion}
			Let $A,B$ be trees spanning a finite subgraph $H\subseteq G$. Let $s,t\in V(H)$ be distinct, with $st\in E(H)$, and let $<$ be a linear order of $V(H)$ with minimum $s$ and maximum $t$. Suppose that $A,B$ and $<$ satisfy conditions~\ref{cond:increasingA} and~\ref{cond:increasingB}. They have extensions $A',B'$ and an extended linear order with the same extrema on a finite subgraph $H'\subseteq G$ satisfying conditions~\ref{cond:increasingA} and~\ref{cond:increasingB} such that, for every component $C$ of $G-V(H')$, the set $N_G(C)$ lies in one component of $A'-E(A)$ and in one component of $B'-E(B)$.
		\end{claim}
		
		\begin{proof}
			In a tree extending $A$, each component remaining after deletion of $E(A)$ contains exactly one vertex of $H$. It contains at least one, since the path from any of its vertices to $H$ uses no edge of $A$ before reaching $H$. It cannot contain two: a path between them avoiding $E(A)$, together with their path in $A$, would contain a cycle. For each vertex $x$ of an extension, let $\alpha(x)\in V(H)$ be this unique vertex. Define $\beta(x)$ in the same way for extensions of $B$. Further tree extensions preserve both labels on vertices already included.
			
			When a path with endvertices $a<b$ is added by~\eqref{eq:ear}, every new vertex has $\alpha$-label $\alpha(a)$ and $\beta$-label $\beta(b)$. We shall use this rule to make both labels constant on the neighbourhood of every remaining component.
			
			There are only finitely many components of $G-V(H)$. Indeed, each has an edge to $H$, whereas only finitely many edges leave the finite vertex set $V(H)$, by local finiteness. We treat these components separately. If a component is finite, all its vertices can be included by finitely many applications of Theorem~\ref{thm:fan} and~\eqref{eq:ear}, with the internal vertices of each added path in that component. To justify this restriction, take a fan from a vertex still omitted to the current vertex set and truncate its paths at their first vertices in that set. Before reaching that set, both paths remain in the component under consideration.
			
			Fix an infinite component $D$ of $G-V(H)$, and let $Z=N_G(V(H))\cap V(D)$.
			The set $Z$ is finite by local finiteness. Since $A$ spans $H$, the graph $H$ is connected; it has at least two vertices by the hypotheses of Claim~\ref{cl:finitecompletion}. Claim~\ref{cl:finiteblocks} therefore applies to $D$. Denote its finite block-cutvertex tree by $\mathcal B$.
			
			Two features guide the block-by-block construction below. First, once a $2$-connected block contains at least two vertices of the current finite subgraph, all further local extensions in that block are obtained by applying Theorem~\ref{thm:fan} inside the block. Their new vertices therefore stay in that block. Second, if the two labels are constant on the included vertices of a block apart from a possible exceptional vertex, we deliberately include all neighbours of that vertex in the block. This is a finite requirement by local finiteness. The exceptional vertex then has no neighbour among the vertices still omitted from the block. These two features will ensure that every remaining component has a neighbourhood on which both labels are constant.
			
			There are edges $az_1,bz_2$ between $H$ and $D$ such that $a\neq b$ and $z_1\neq z_2$. Indeed, if no two edges in this bipartite graph had disjoint ends, all its edges would share one endvertex. To see this, take an edge $az$; if there are edges $a'z$ with $a'\neq a$ and $az'$ with $z'\neq z$, these two edges have disjoint ends. Otherwise all edges share $a$ or all share $z$. Either possibility would give a cutvertex of $G$: the common endvertex would separate $D$ from $H$, or the nonempty graph $D-z$ from $H$. Here $H$ has at least two vertices, as required by the ordering conditions.
			
			Join $z_1$ to $z_2$ by a path in $D$. Together with $az_1$ and $bz_2$, it gives an $H$-path $P$. Reverse its notation if necessary so that $a<b$, and add it using~\eqref{eq:ear}. All its internal vertices then have label pair $(a,b)$ and occur increasingly from $a$ to $b$.
			
			The blocks containing edges of the $z_1$--$z_2$ subpath occur along a path in $\mathcal B$. In each such block $K$, this subpath has a nontrivial connected intersection. Indeed, leaving and later re-entering a block would require either repeating a cutvertex on the graph path or making a cycle in $\mathcal B$, contrary to Theorem~\ref{thm:blocks}. If $K$ is a bridge, both its vertices are already present. Otherwise $K$ is $2$-connected. Theorem~\ref{thm:fan}, applied inside $K$, permits any specified omitted vertex of $K$ to be included by a path whose distinct ends are among the vertices of $K$ already present and whose internal vertices are new. There are at least two vertices of $K$ already present, because the original intersection is nontrivial. Apply~\eqref{eq:ear} to each such path.
			
			In this way, include all cutvertices of $D$ in $K$ and all vertices of $Z\cap V(K)$. These prescribed sets are finite. Every vertex added inside $K$ still has label pair $(a,b)$, because the two endvertices of every local extension have that pair. Carry out these finite extensions in all the blocks traversed by $P$, and mark these blocks as processed.
			
			At each completed step, let $W$ be the current vertex set. We maintain the following two properties of the processed blocks:
			\begin{enumerate}[label=(B\arabic*),ref=(B\arabic*),leftmargin=*]
				\item\label{cond:processedblocks} The processed blocks, together with all cutvertices of $D$ belonging to them, induce a connected subtree of $\mathcal B$. Every vertex of $W\cap V(D)$ belongs to a processed block. For every processed block $J$, the set $W$ contains all cutvertices of $D$ in $J$ and all vertices of $Z\cap V(J)$.
				\item\label{cond:blocklabels} For every processed block $J$, the two labels are constant on $W\cap V(J)$, with at most one exceptional vertex $v_J$. If such a vertex is present, it belongs to $W$ and all its neighbours in $J$ belong to $W$; the labels are constant on $(W\cap V(J))\setminus\{v_J\}$.
			\end{enumerate}
			The extensions made so far satisfy condition~\ref{cond:processedblocks} by construction and condition~\ref{cond:blocklabels} with no exceptional vertex. Every unprocessed block meets $W$ in at most one vertex. Indeed, any vertex in that intersection belongs to a processed block by condition~\ref{cond:processedblocks}; it is therefore a cutvertex by Theorem~\ref{thm:blocks}. Two such vertices would give two connections between the unprocessed block and the connected processed subtree, creating a cycle in $\mathcal B$.
			
			Suppose that a block remains unprocessed. Choose an unprocessed block $K$ incident with the processed subtree at a cutvertex $v$, and let $U$ be the component of $D-v$ containing $K-v$. The separation statement in Theorem~\ref{thm:blocks} shows that $U$ contains no vertex of the current subgraph. Since $G-v$ is connected, $U$ has a neighbour in the original $H$: its only possible edges to $D-U$ are incident with $v$, and all its neighbours outside $D$ lie in $H$. Choose a path from $v$ into $K-v$, through $U$, and then to $H$, stopping at its first vertex in $H$. It is a path of the form used in~\eqref{eq:ear}. Within $D$ it traverses a path of unprocessed blocks starting with $K$. Add it by~\eqref{eq:ear}; all its new vertices receive one common pair of labels.
			
			For every newly traversed block other than $K$, the vertices now present are exactly those of its nontrivial intersection with the new path. Apply the same finite extensions inside each of these blocks to include all cutvertices of $D$ in that block and all its vertices in $Z$. All its included vertices have the common pair of labels.
			
			In $K$, the cutvertex $v$ of $D$ may have a different label from the other vertices on the new path. Write $Q$ for the part of that path in $K$. If $v$ is the smaller endvertex of the whole extension, it is the minimum vertex of $Q$. Its $\alpha$-label agrees with that of $Q-v$, whereas its $\beta$-label may differ. If $v$ is the larger endvertex, it is the maximum vertex of $Q$, and only its $\alpha$-label may differ.
			
			Extend the trees inside $K$ to include all cutvertices of $D$ in $K$, all vertices of $Z\cap V(K)$, and all neighbours of $v$ in $K$. This last prescribed set is finite by local finiteness. At every such extension the new vertices are inserted between its endvertices, so all vertices included in $K$ remain between the minimum and maximum of $Q$. If $v$ is the minimum, it can never be the larger endvertex of a local extension. The second tree therefore attaches each new path at a vertex other than $v$, preserving the common $\beta$-label on $K-v$; the $\alpha$-label is already constant throughout the included part of $K$. If $v$ is the maximum, it can never be the smaller endvertex, and the same argument with the two trees interchanged applies. Hence all included vertices of $K-v$ have the common pair of labels. Also, no vertex of $K$ still omitted is adjacent to $v$, because all neighbours of $v$ in $K$ have been included. For a bridge $K$, these conclusions hold immediately, since both endvertices are present.
			
			Figure~\ref{fig:blockcompletion} shows the block traversal and the treatment of the possible exceptional vertex $v$.
			
			\begin{figure}[htbp]
				\centering
				\begin{tikzpicture}[
					x=1cm,y=1cm,
					every node/.style={font=\fontsize{10}{12}\selectfont},
					vtx/.style={circle,fill=black,inner sep=1.4pt},
					selected/.style={draw=black,line width=.9pt},
					omitted/.style={draw=gray!60,dashed,line width=.45pt},
					oldregion/.style={draw=gray!50,fill=gray!6,rounded corners=5pt},
					block/.style={rectangle,draw=black,fill=white,minimum width=.7cm,minimum height=.45cm,inner sep=2pt}
					]
					\node at (2.7,3.4) {(a) Blocks of $D$};
					\node[block] (B0) at (.5,.7) {$B_0$};
					\node[vtx,label=above:$x$] (x) at (1.5,.7) {};
					\node[block] (B1) at (2.5,.7) {$B_1$};
					\node[vtx,label=left:$v$] (v) at (2.5,1.55) {};
					\node[block] (K) at (2.5,2.45) {$K$};
					\node[vtx,label=above:$w$] (w) at (3.7,2.45) {};
					\node[block] (L) at (4.9,2.45) {$L$};
					\begin{scope}[on background layer]
						\node[oldregion,fit=(B0)(x)(B1)(v),inner sep=6pt] {};
					\end{scope}
					\draw[selected] (B0)--(x)--(B1);
					\draw[draw=gray!70,line width=.45pt] (B1)--(v);
					\draw[omitted] (v)--(K)--(w)--(L);
					\node[align=center,anchor=north] at (1.5,-.08) {blocks traversed by\\the initial $H$-path};
					\node[align=center] at (4.25,1.45) {next blocks\\traversed};
					
					\begin{scope}[xshift=6.2cm]
						\node at (2.2,3.4) {(b) The block $K$};
						\draw[oldregion,fill=gray!3] (-.35,.25) rectangle (4.8,2.95);
						\node[anchor=north west] at (-.25,2.9) {$K$};
						\draw[rounded corners=8pt,fill=gray!12,draw=gray!55]
						(2.45,.65) rectangle (4.55,2.65);
						\node at (4.25,2.4) {$C$};
						\node[vtx,label=below:$v$] (vdetail) at (0,1.6) {};
						\node[vtx,label=above:$p$] (p) at (1.2,2.35) {};
						\node[vtx,label=below:$q$] (q) at (1.2,.85) {};
						\node[vtx,label=above:$w$] (wdetail) at (1.75,1.6) {};
						\coordinate (c1) at (2.85,2.05);
						\coordinate (c2) at (2.85,1.05);
						\coordinate (c3) at (3.85,1.55);
						\draw[selected] (vdetail)--(p)--(wdetail)--(q)--(vdetail)--(wdetail);
						\draw[omitted] (p)--(c1)--(c3)--(c2)--(q);
						\draw[omitted] (wdetail)--(c1) (wdetail)--(c2);
						\foreach \c in {c1,c2,c3}{
							\node[circle,draw=gray!65,fill=white,inner sep=1.4pt] at (\c) {};
						}
						\node[align=center,anchor=north] at (2.2,-.08)
						{$N_G(C)\subseteq\{p,q,w\}$\\$p,q,w$ have the same two labels};
					\end{scope}
				\end{tikzpicture}
				\caption{The construction in Claim~\ref{cl:finitecompletion}. In (a), the shaded subtree is processed, the boxes are blocks and the dots are cutvertices. The vertex $w$ belongs to both $K$ and $L$. In (b), all neighbours of $v$ in $K$ are included, leaving $N_G(C)\subseteq\{p,q,w\}$ with constant labels. Black edges in (b) indicate adjacency among the included vertices.}
				\label{fig:blockcompletion}
			\end{figure}
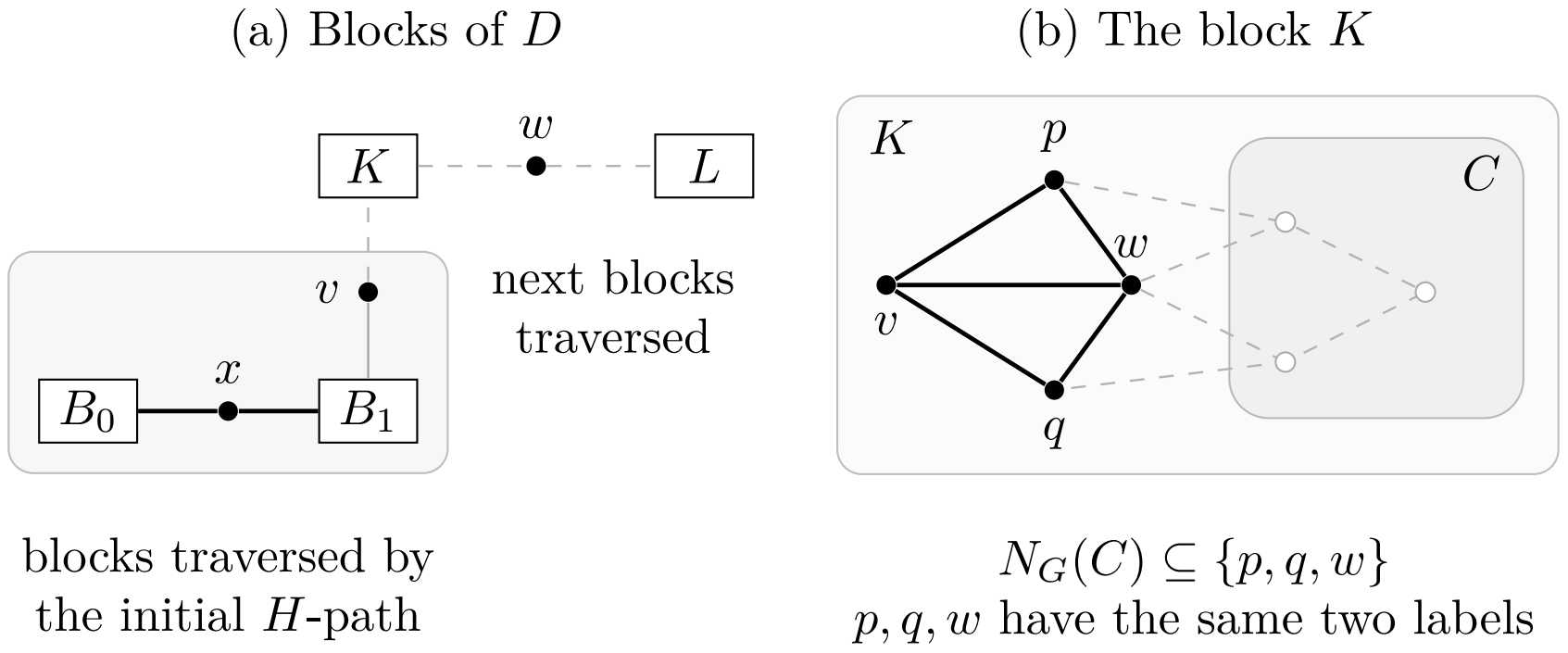
			
			Mark all blocks traversed at this step as processed. The newly processed blocks and the cutvertices of $D$ in them, together with the vertices of the previously processed subtree, induce a connected subtree of $\mathcal B$. The extensions inside these blocks have included all cutvertices of $D$ belonging to them and all their vertices in $Z$, so condition~\ref{cond:processedblocks} is preserved. Extending the trees inside one newly traversed block cannot introduce an additional vertex into another block on the same path: their only possible common vertex is the cutvertex already used by that path, by Theorem~\ref{thm:blocks}.
			
			Nor can any of these extensions introduce a new vertex into an earlier processed block. An intersection with such a block is a cutvertex already included by condition~\ref{cond:processedblocks}; a second connection would create a cycle in $\mathcal B$. Thus the included vertex sets and their labels in earlier processed blocks do not change. They continue to satisfy condition~\ref{cond:blocklabels}. Each newly traversed block other than $K$ has constant labels on its included vertices, while the local construction in $K$ permits only the exception $v$ and includes every neighbour of $v$ in $K$. This proves condition~\ref{cond:blocklabels} for all newly processed blocks as well. Every unprocessed block again meets the current vertex set in at most one vertex by condition~\ref{cond:processedblocks} and Theorem~\ref{thm:blocks}. Each step processes at least one block, and $\mathcal B$ is finite by Claim~\ref{cl:finiteblocks}, so the procedure terminates after finitely many finite extensions.
			
			Let $W$ be the final vertex set. Every block of $D$ is now processed, so condition~\ref{cond:processedblocks} implies that all cutvertices of $D$ and all vertices of $Z$ belong to $W$. Let $C$ be a component of $D-W$. Each vertex of $C$ belongs to a unique block of $D$, since all cutvertices lie in $W$. By Theorem~\ref{thm:blocks}, every edge of $D$ belongs to a block. It follows by connectedness that the vertices of $C$ lie in a single block $K$, and every edge of $D$ incident with $C$ belongs to $K$. Thus $N_D(C)\subseteq W\cap V(K)$. Also, $C$ has no neighbour in the original $H$, since $Z\subseteq W$. By condition~\ref{cond:blocklabels}, both labels are constant on $W\cap V(K)$ apart from at most one vertex. If the exceptional vertex $v_K$ exists, all its neighbours in $K$ belong to $W$, so $v_K$ has no neighbour in $C$. It follows that both labels are constant on $N_G(C)$, as required.
			
			Constructions in different components of $G-V(H)$ have disjoint internal vertex sets and use only vertices of the original $H$ outside their respective components. They therefore preserve the conclusions obtained in earlier components. There are finitely many components and all extensions are finite. Their union gives $H',A',B'$ as required by Claim~\ref{cl:finitecompletion}.
		\end{proof}
		
		\begin{claim}\label{cl:finitecuts}
			Let $H,A,B$ and $H',A',B'$ be as in Claim~\ref{cl:finitecompletion}. In every spanning tree $T$ of $G$ containing $A'$, the fundamental cut of each edge of $A$ is finite. The corresponding assertion holds for $B$ and every spanning tree containing $B'$.
		\end{claim}
		
		\begin{proof}
			It suffices to prove the assertion for $A$. Fix a spanning tree $T$ containing $A'$. Label each vertex of $H'$ by the unique vertex of $H$ in its component of $A'-E(A)$, as in the proof of Claim~\ref{cl:finitecompletion}. For each component $C$ of $G-V(H')$, give every vertex of $C$ the common label of $N_G(C)$ provided by Claim~\ref{cl:finitecompletion}.
			
			First observe that $T[V(H')]=A'$. Indeed, any additional edge of $T$ with both ends in $H'$ would form a cycle together with the path between its ends in $A'$. Hence every edge of $T$ with both ends in $H'$ and differently labelled ends belongs to $A$: all edges of $A'-E(A)$ lie within the components used to define the labels. Every edge with at least one end outside $H'$ has equally labelled ends, by the assignment of a common label to each component $C$ and its neighbourhood. Consequently an edge of $T$ whose ends have different labels must be an edge of $A$.
			
			The label classes are exactly the components of $T-E(A)$. The preceding observation shows that no component can contain two different labels. Conversely, the vertices of $H'$ with any fixed label are connected in $A'-E(A)$. For a vertex outside $H'$, let $U$ be its component in $T-V(H')$. Since $T$ is connected, $N_T(U)\cap V(H')\neq\varnothing$. This set contains exactly one vertex, say $u$: a path through $U$ joining two distinct vertices of $H'$, together with their path in $A'$, would create a cycle in $T$. Let $C$ be the component of $G-V(H')$ containing $U$. Then $u\in N_G(C)$, so $u$ has the label assigned to every vertex of $U$. The subgraph of $T$ induced by $V(U)\cup\{u\}$ is connected and contains no edge of $A$. Thus every vertex of $U$ belongs to the component of $T-E(A)$ containing the vertices of $H'$ with its label, proving the assertion.
			
			Every edge of $G$ whose ends have different labels has both ends in $H'$. An edge with both ends outside $H'$ lies within one component of $G-V(H')$, where the label is constant. An edge with exactly one end in $H'$ joins that end to a component $C$ whose label agrees with every vertex of $N_G(C)$. Such an edge also has equally labelled ends.
			
			Now let $e\in E(A)$. The two components of $T-e$ are unions of components of $T-E(A)$. Thus every edge in the fundamental cut of $e$ joins vertices with different labels and consequently has both ends in $H'$. This cut is contained in the edge set of the finite graph $G[V(H')]$, and is therefore finite. The proof for $B$ is identical.
		\end{proof}
		
		Figure~\ref{fig:fundamentalcut} illustrates why the remaining components contribute no further edge to the fundamental cut.
		
		\begin{figure}[htbp]
			\centering
			\begin{tikzpicture}[
				x=1cm,y=1cm,
				every node/.style={font=\fontsize{10}{12}\selectfont},
				vtx/.style={circle,fill=black,inner sep=1.4pt},
				selected/.style={draw=black,line width=.9pt},
				omitted/.style={draw=gray!60,dashed,line width=.45pt},
				oldregion/.style={draw=gray!50,fill=gray!6,rounded corners=5pt}
				]
				\draw[oldregion] (-5.4,0) rectangle (5.4,1.6);
				\node[anchor=north west] at (-5.3,1.55) {$H'$};
				\draw[oldregion] (-4.7,-2.15) rectangle (-1.9,-.9);
				\draw[oldregion] (1.9,-2.15) rectangle (4.7,-.9);
				\node at (-3.3,-1.93) {$C_1$};
				\node at (3.3,-1.93) {$C_2$};
				\draw[densely dotted,line width=.65pt] (0,-2.4)--(0,2.1);
				\node[align=center] at (0,2.48) {partition defined by $T-e$};
				\node at (-2.9,1.93) {$\alpha=a$};
				\node at (2.9,1.93) {$\alpha=b$};
				\node at (-3.3,-2.48) {$\alpha=a$};
				\node at (3.3,-2.48) {$\alpha=b$};
				\node[vtx,label=below:$a$] (a) at (-2.8,.5) {};
				\node[vtx,label=below:$b$] (b) at (2.8,.5) {};
				\node[vtx] (a1) at (-4.6,.5) {};
				\node[vtx] (a2) at (-3.5,1.3) {};
				\node[vtx] (b1) at (4.6,.5) {};
				\node[vtx] (b2) at (3.5,1.3) {};
				\node[vtx] (c11) at (-4,-1.35) {};
				\node[vtx] (c12) at (-2.7,-1.55) {};
				\node[vtx] (c21) at (4,-1.35) {};
				\node[vtx] (c22) at (2.7,-1.55) {};
				\draw[selected] (a1)--(a)--(a2) (b1)--(b)--(b2);
				\draw[selected] (a)--(b);
				\node at (-.4,.85) {$e$};
				\draw[selected] (a1)--(c11)--(c12) (b1)--(c21)--(c22);
				\draw[omitted] (a2)--(b2);
				\draw[omitted] (a2) to[bend left=12] (c12);
				\draw[omitted] (b2) to[bend right=12] (c22);
			\end{tikzpicture}
			\caption{Finite fundamental cuts in Claim~\ref{cl:finitecuts}. Black edges belong to $T$; dashed edges belong only to $G$. The label is constant on $N_G(C)$ for each component $C$ of $G-V(H')$. Thus every edge in the fundamental cut has both ends in $H'$. Only portions of the components of $T-E(A)$ labelled $a,b$ are shown.}
			\label{fig:fundamentalcut}
		\end{figure}
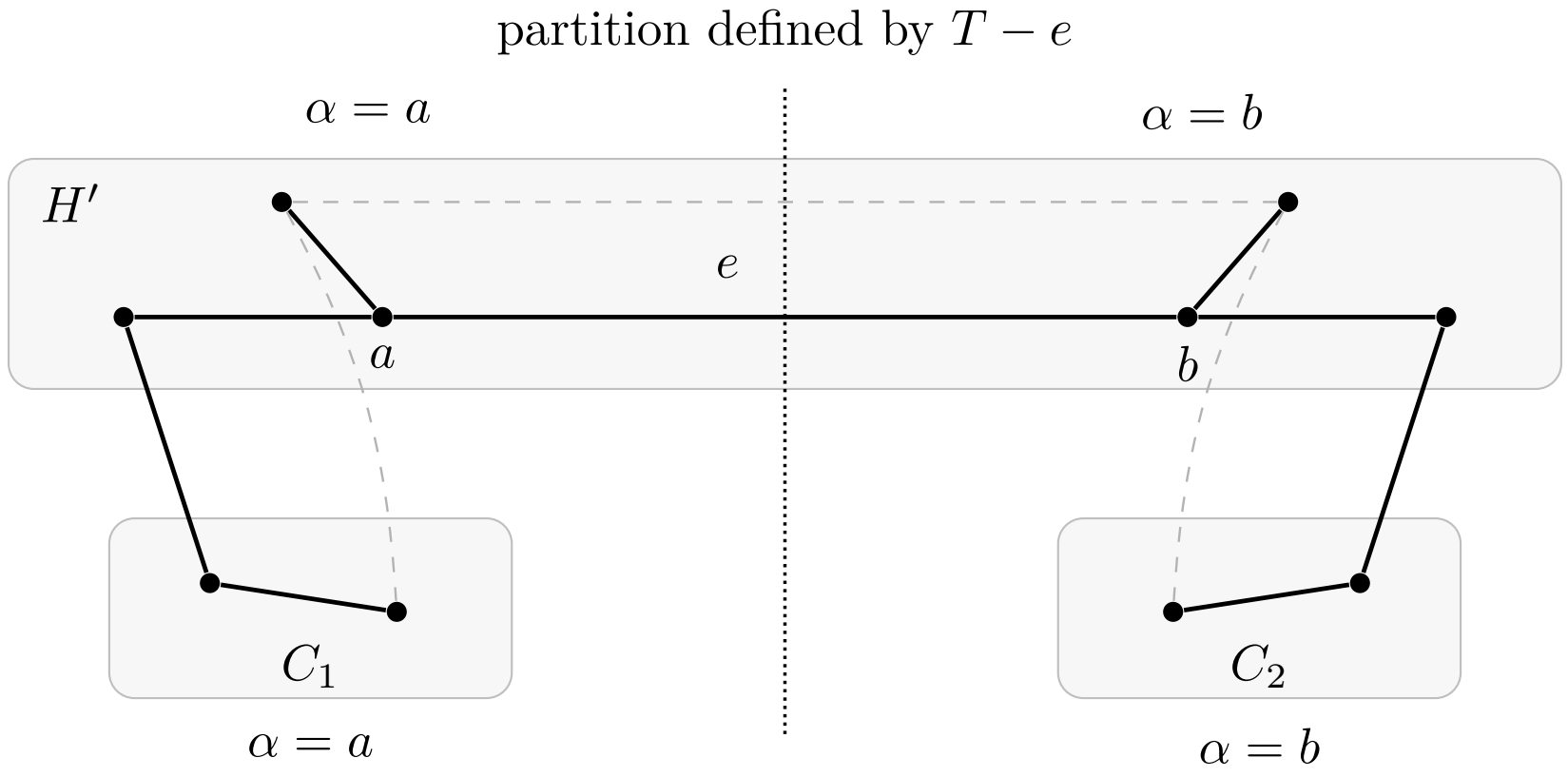
		
		Fix the prescribed root $s$. Since $G$ is locally finite and connected, it is countable. Enumerate its vertices as $V(G)=\{v_0,v_1,\ldots\}$. Choose a cycle $s u_1\ldots u_m t s$, as in the proof of Theorem~\ref{thm:two}, and start with the trees obtained by deleting $st$ and $su_1$, respectively, and the order $s<u_1<\cdots<u_m<t$.
		
		At stage $n$, include $v_n$ if it is still omitted, using Theorem~\ref{thm:fan} and~\eqref{eq:ear}, and then apply Claim~\ref{cl:finitecompletion}. Each stage adds finitely many vertices and preserves conditions~\ref{cond:increasingA} and~\ref{cond:increasingB}. Take the unions of the first and second trees over all stages, denoting them by $T_1,T_2$. Every vertex appears at some stage and retains its path to $s$. Each union is therefore connected and spanning; a cycle in a union would already occur at a finite stage, which is impossible. Thus the unions are trees, and the orders extend to a linear order of their common vertex set. For $v\notin\{s,t\}$, conditions~\ref{cond:increasingA} and~\ref{cond:increasingB} place all internal vertices of $sT_1v$ below $v$ and all internal vertices of $sT_2v$ above $v$. These paths are internally disjoint and distinct, since $s$ is a leaf of $T_2$ adjacent to $t$. For the target $t$, the second path is $st$, whereas the first is the preserved initial path $s u_1\ldots u_m t$. They too are distinct and internally disjoint. Thus $T_1,T_2$ are independent.
		
		Every edge of $T_i$ is present before some application of Claim~\ref{cl:finitecompletion}; since the final tree contains the resulting finite extension, Claim~\ref{cl:finitecuts} shows that its fundamental cut in $T_i$ is finite. Both trees are therefore end-faithful by Proposition~\ref{prop:finitecutcriterion}.
	\end{proof}
	
	The local finiteness assumption in Theorem~\ref{thm:faithful} is essential, even when only one vertex has infinite degree.
	
	\begin{example}\label{ex:fan}
		Let $G$ be obtained from a ray $R=v_0v_1\ldots$ by adding a vertex $z$ adjacent to every vertex of $R$. The graph is shown in Figure~\ref{fig:fan}. Notice that $G$ is $2$-connected and has a unique end, but has no two independent end-faithful spanning trees rooted at any vertex. Only $z$ has infinite degree.

		\begin{figure}[htbp]
			\centering
			\begin{tikzpicture}[
				x=1.25cm,y=1cm,
				every node/.style={font=\fontsize{10}{12}\selectfont},
				vtx/.style={circle,fill=black,inner sep=1.4pt},
				selected/.style={draw=black,line width=.9pt},
				omitted/.style={draw=gray!60,dashed,line width=.45pt}
				]
				\node[vtx,label=above:$z$] (zfan) at (2.5,2.25) {};
				\foreach \i in {0,...,5}{
					\node[vtx,label=below:$v_{\i}$] (fan\i) at (\i,0) {};
					\draw[draw=gray!65,line width=.55pt] (zfan)--(fan\i);
				}
				\foreach \i [evaluate=\i as \j using int(\i+1)] in {0,...,4}{
					\draw[selected] (fan\i)--(fan\j);
				}
				\draw[selected] (fan5)--(5.55,0);
				\node at (5.85,0) {$\cdots$};
				\draw[omitted] (zfan)--(5.85,.25);
			\end{tikzpicture}
			\caption{The graph of Example~\ref{ex:fan} is a ray $R=v_0v_1\ldots$ with a universal vertex $z$. Any two rays have a common tail in $R$, so they have infinitely many vertices in common. This rules out two independent end-faithful spanning trees.}
			\label{fig:fan}
		\end{figure}
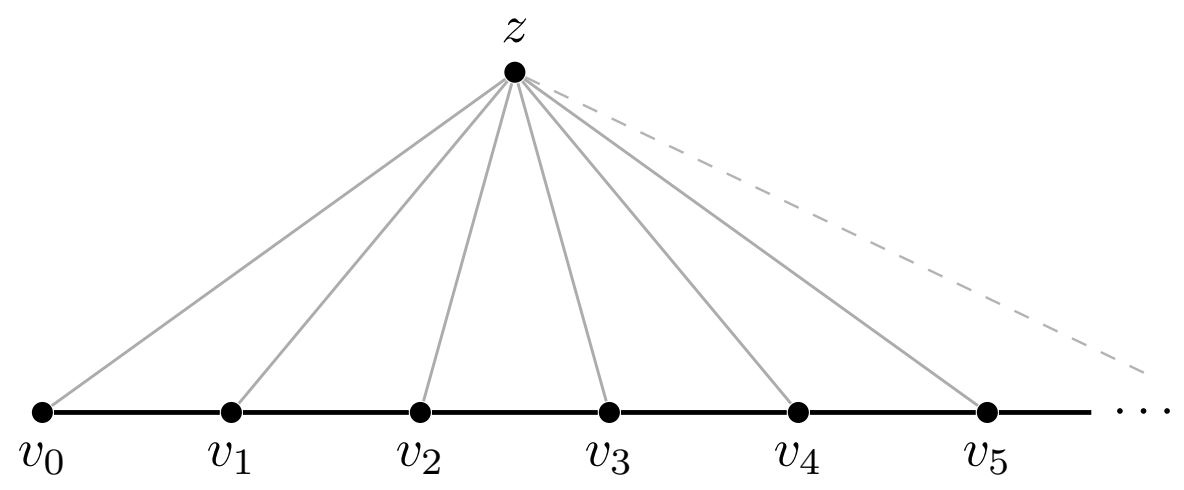
		
		Deleting $z$ leaves $R$ connected, and deleting any $v_i$ leaves the remaining vertices connected through $z$. Thus $G$ is $2$-connected. Deleting $z$ and $v_1$ isolates $v_0$, so its vertex-connectivity is exactly two. Every ray of $G$ has a tail in $R$, since it meets $z$ at most once. A ray contained in $R$ is a tail of $R$. Hence $G$ has a unique end, and any two rays of $G$ have infinitely many vertices in common.
		
		If two independent trees have the same root, any two rays, one in each tree and both starting at that root, have only finitely many common vertices. Indeed, let $P_i\subseteq T_i$ and $P_j\subseteq T_j$ be rays starting at the common root $r$. If common vertices $x,y\neq r$ occur in this order on $P_i$, independence of the root paths to $y$ forces them to occur in the opposite order on $P_j$. Otherwise $x$ would be an internal vertex of both root paths to $y$. Infinitely many common vertices would give an infinite strictly decreasing sequence of their positions on $P_j$, which is impossible. In particular, any finite family of rays starting at the common root in distinct independent trees has pairwise disjoint tails: for each ray, delete a finite initial segment containing all its intersections with the other rays.
		
		If two independent spanning trees of $G$ were both end-faithful, each would contain a ray starting at their common root. This contradicts the preceding observation, since any two rays of $G$ have infinitely many common vertices. The argument applies to every choice of root.
	\end{example}
	
	\section{Three end-faithful spanning trees in chordal graphs}\label{sec:chordal}
	
	We now prove  that every locally finite $3$-connected chordal graph has three independent end-faithful spanning trees rooted at any prescribed vertex.
	
	\begin{definition}\label{def:treedecomposition}
		A \emph{tree-decomposition} of a graph $G$ is a pair $(\mathcal D,(B_t)_{t\in V(\mathcal D)})$, where $\mathcal D$ is a tree and each $B_t$ is a subset of $V(G)$. The sets $B_t$ cover $V(G)$, every edge of $G$ has both endvertices in some $B_t$, and for every $v\in V(G)$ the set $\{t\in V(\mathcal D):v\in B_t\}$ induces a connected subtree of $\mathcal D$. The sets $B_t$ are the \emph{bags} of the decomposition. For $tu\in E(\mathcal D)$, the intersection $B_t\cap B_u$ is the \emph{adhesion set} associated with $tu$. A \emph{clique tree} of $G$ is a tree-decomposition whose bags are exactly the maximal cliques of $G$, each occurring once.
	\end{definition}

	We use the following property of tree-decompositions.
	
	\begin{proposition}[Diestel~{\cite[Lemma~12.3.1]{Diestel}}]\label{prop:decompositionseparation}
		Let $(\mathcal D,(B_t)_{t\in V(\mathcal D)})$ be a tree-decomposition of a graph $G$. Given $tu\in E(\mathcal D)$, let $\mathcal D_t$ and $\mathcal D_u$ be the components of $\mathcal D-tu$ containing $t$ and $u$, respectively. For $i\in\{t,u\}$, let $X_i$ be the union of the bags indexed by $V(\mathcal D_i)$.
		Then $X_t\cup X_u=V(G)$, $X_t\cap X_u=B_t\cap B_u$, and no edge of $G$ joins $X_t\setminus X_u$ to $X_u\setminus X_t$. In particular, every path from $X_t\setminus X_u$ to $X_u\setminus X_t$ meets $B_t\cap B_u$.
	\end{proposition}

	We use the following result of Hofer-Temmel and Lehner.
	
	\begin{theorem}[Hofer-Temmel and Lehner~{\cite[Theorem~3]{HTL}}]\label{thm:cliquetree}
		Every connected locally finite chordal graph has a clique tree.
	\end{theorem}

	\begin{theorem}\label{thm:chordal}
		Every locally finite $3$-connected chordal graph has three independent end-faithful spanning trees rooted at any prescribed vertex.
	\end{theorem}
	
	\begin{proof}
		If $G$ is finite, apply Theorem~\ref{thm:finite} with $k=3$; end-faithfulness is vacuous. Suppose that $G$ is infinite. Fix $r\in V(G)$ and choose a clique tree by Theorem~\ref{thm:cliquetree}. Since each maximal clique occurs once, we may identify a vertex of its decomposition tree $\mathcal D$ with its bag. Root $\mathcal D$ at a bag $B_0$ containing $r$.
		
		Every bag is finite. Moreover, each vertex $v$ of $G$ belongs to only finitely many bags, since any clique containing $v$ is a subset of the finite set $\{v\}\cup N_G(v)$. For adjacent bags $B,P$, Proposition~\ref{prop:decompositionseparation} shows that their adhesion set $S=B\cap P$ separates two nonempty sets: these contain $B\setminus P$ and $P\setminus B$, respectively, and both sets are nonempty by maximality. Hence $|S|\geq3$, by the $3$-connectivity of $G$. Each of $B,P$ properly contains $S$, so every bag has at least four vertices. The tree $\mathcal D$ has more than one vertex, since a single finite bag cannot cover the infinite graph $G$.
		
		The tree $\mathcal D$ is locally finite. Indeed, every neighbour of a bag $B$ meets $B$, and $B$ meets only finitely many bags: it is finite, and each of its vertices belongs to finitely many bags. In particular, every level of $\mathcal D$ is finite. Process its bags in nondecreasing distance from $B_0$, with an arbitrary order within each level. Every bag is then reached after finitely many steps. If $B\neq B_0$ has parent $P$, then
		\begin{equation}\label{eq:cliqueattachment}
			B\cap\bigcup_{A\text{ preceding }B}A=B\cap P.
		\end{equation}
		This is the running intersection property for the chosen enumeration of the clique tree~\cite[Theorem~18 and the following paragraph]{HTL}. Thus the vertices of $B$ already present when it is processed are those in $S=B\cap P$, and $B\setminus S\neq\varnothing$.
		
		We construct three independent trees $T_1,T_2,T_3$ on the union of the processed bags, always adding new vertices as leaves. For each processed bag $B$, we choose distinct vertices $b_1(B),b_2(B),b_3(B)\in B$ and maintain the following conditions:
		\begin{enumerate}[label=(C\arabic*),ref=(C\arabic*),leftmargin=*]
			\item\label{cond:bagancestor} For each $i\in\{1,2,3\}$ and $x,y\in B$, if $x\neq r$ is a proper ancestor of $y$ in $T_i$, then $x=b_i(B)$.
			\item\label{cond:bagpaths} For distinct $i,j\in\{1,2,3\}$ and all $x,y\in B$, $V(rT_ix)\cap V(rT_jy)\subseteq B\cup\{r\}$.
		\end{enumerate}
		The chosen vertices may include $r$. Once these conditions hold for a processed bag, they remain true when leaves are added, since root paths to old vertices do not change.
		
		For each $i$, we also define a linear order $\prec_i$ on the vertices of $T_i$. Give the children of each vertex distinct positive integer indices, assigning a new child the next unused index. Associate each vertex with the finite sequence of indices along its root path; the sequence for $r$ is empty. Compare these sequences lexicographically: at the first position where they differ, the sequence with the smaller entry comes first; if one is a proper prefix of the other, the shorter sequence comes first. Distinct vertices have distinct sequences, so this is a linear order. Adding leaves changes neither an old sequence nor the comparison of old vertices. Each vertex precedes its proper descendants. If $x\prec_i y$ and $x,y$ have the same depth, every descendant of $x$ precedes every descendant of $y$: the first differing entries in the sequences for $x$ and $y$ also determine the comparison of any extensions of those sequences.
		
		Choose distinct $a_1,a_2,a_3\in B_0\setminus\{r\}$. Initially, let $T_i$ be the star on $B_0$ with centre $a_i$, rooted at $r$, and set $b_i(B_0)=a_i$. Assign indices to its children arbitrarily. For $v\neq r$, the path $rT_iv$ is $ra_i$ if $v=a_i$, and $ra_iv$ otherwise. Thus the three root paths to $v$ are distinct and internally disjoint. If a vertex different from $r$ is a proper ancestor in $T_i$ of a vertex of $B_0$, it must be $a_i$, so condition~\ref{cond:bagancestor} holds. Condition~\ref{cond:bagpaths} holds because all these paths are contained in $B_0$.
		
		Choosing the least available vertex below will ensure that, for each fixed depth $d$, all vertices in sufficiently distant bags along any ray of $\mathcal D$ from $B_0$ have the same ancestor at depth $d$ in $T_i$; see Claim~\ref{cl:descendants}.
		
		Now suppose that $B$ is the next bag, with parent $P$, and write $S=B\cap P$. Choose distinct $c_1,c_2,c_3\in S$ as follows. First, simultaneously for every $i$ such that $b_i(P)\in S$, set $c_i=b_i(P)$. Then consider the remaining indices in the order $1,2,3$, choosing $c_i$ to be the $\prec_i$-least vertex of $S$ not already selected. The first choices are distinct by induction, and every subsequent choice is possible because $|S|\geq3$.
		
		The paths $rT_1c_1,rT_2c_2,rT_3c_3$ meet pairwise only at $r$. Indeed, suppose that $w\neq r$ belongs to the paths for distinct indices $i,j$. By condition~\ref{cond:bagpaths} on $P$, we have $w\in P$. If $w\notin\{c_i,c_j\}$, then $w$ is a proper ancestor of $c_i$ in $T_i$ and of $c_j$ in $T_j$. Condition~\ref{cond:bagancestor} on $P$ gives $w=b_i(P)=b_j(P)$, contrary to the distinctness of the chosen vertices. Otherwise we may assume that $w=c_i$. Since $c_i\neq c_j$, condition~\ref{cond:bagancestor} in $T_j$ gives $w=b_j(P)\in S$. But the first step of the choice then gives $c_j=b_j(P)=w=c_i$, again a contradiction.
		
		For every $v\in B\setminus S$ and each $i$, add the edge $c_iv$ to $T_i$, and set $b_i(B)=c_i$. These edges exist because $B$ is a clique. Assign indices to the new children at $c_i$. Each new vertex is a leaf, so the trees remain trees and their old root paths remain unchanged. The root path to a new vertex $v$ consists of $rT_ic_i$ followed by $c_iv$. The preceding paragraph and the distinctness of the $c_i$ show that these three paths are distinct and internally disjoint, including when one $c_i$ equals $r$. Thus independence is preserved.
		
		We verify condition~\ref{cond:bagancestor} on $B$. Fix $i$ and suppose that $x\in B\setminus\{r\}$ is a proper ancestor of $y\in B$ in $T_i$. If $y\in S$, its root path is old, so $x\in S$ by Equation~\ref{eq:cliqueattachment}. Condition~\ref{cond:bagancestor} on $P$ gives $x=b_i(P)\in S$, which was retained as $c_i=b_i(B)$. If $y\in B\setminus S$, then either $x=c_i$, as required, or $x$ is a proper ancestor of $c_i$ on its old root path. In the latter case $x\in S$, again by equation~\ref{eq:cliqueattachment}, and condition~\ref{cond:bagancestor} on $P$ gives $x=b_i(P)\in S$. The choice of $c_i$ would imply $c_i=x$, contradicting the assumption that $x$ is a proper ancestor of $c_i$.
		
		For condition~\ref{cond:bagpaths}, fix distinct $i,j$ and $x,y\in B$. Put $x'=x$ if $x\in S$ and $x'=c_i$ otherwise; likewise, put $y'=y$ if $y\in S$ and $y'=c_j$ otherwise. Thus $x',y'\in S$. Any common vertex $w\notin B\cup\{r\}$ of $rT_ix$ and $rT_jy$ lies on the old paths $rT_ix'$ and $rT_jy'$. Condition~\ref{cond:bagpaths} on $P$ implies $w\in P\setminus S$. Since $x',y'\in S$, the vertex $w\neq r$ is a proper ancestor of $x'$ in $T_i$ and of $y'$ in $T_j$. Condition~\ref{cond:bagancestor} gives $w=b_i(P)=b_j(P)$, a contradiction. This completes the inductive step.
		
		The attachment at one new bag is illustrated in Figure~\ref{fig:cliqueextension}.
		
		\begin{figure}[htbp]
			\centering
			\begin{tikzpicture}[
				x=1cm,y=1cm,
				font=\fontsize{10}{12}\selectfont,
				vtx/.style={circle,fill=black,inner sep=1.4pt},
				selected/.style={line width=.9pt},
				omitted/.style={draw=gray!60,dashed,line width=.45pt},
				oldregion/.style={fill=gray!6,rounded corners=4pt},
				first tree/.style={draw=black,line width=.9pt},
				second tree/.style={draw=black,line width=.9pt,dash pattern=on 4pt off 2.5pt},
				third tree/.style={draw=black,line width=.9pt,dash pattern=on .7pt off 2pt},
				every node/.append style={outer sep=1pt}
				]
				\path[oldregion] (3.2,.3) rectangle (7.7,3.8);
				\path[oldregion] (6.1,0) rectangle (11.7,4.1);
				\path[fill=gray!14] (6.1,.3) rectangle (7.7,3.8);
				\draw[gray!55,line width=.45pt,rounded corners=4pt]
				(3.2,.3) rectangle (7.7,3.8);
				\draw[gray!55,line width=.45pt,rounded corners=4pt]
				(6.1,0) rectangle (11.7,4.1);
				\node at (4.1,3.5) {$P$};
				\node at (10.8,3.8) {$B$};
				\node at (6.9,3.5) {$S$};
				\coordinate (r) at (.35,2);
				\coordinate (c1) at (6.9,2.95);
				\coordinate (c2) at (6.9,2);
				\coordinate (c3) at (6.9,1.05);
				\coordinate (v) at (10.5,2);
				\draw[first tree] (r) .. controls (1.8,3.05) and (4,2.95) .. (c1);
				\draw[second tree] (r)--(c2);
				\draw[third tree] (r) .. controls (1.8,.95) and (4,1.05) .. (c3);
				\node at (2.3,3.3) {$rT_1c_1$};
				\node at (2.3,1.65) {$rT_2c_2$};
				\node at (2.3,.65) {$rT_3c_3$};
				\draw[first tree] (c1)--(v);
				\draw[second tree] (c2)--(v);
				\draw[third tree] (c3)--(v);
				\node[vtx,label=left:$r$] at (r) {};
				\node[vtx,label=above left:$c_1$] at (c1) {};
				\node[vtx,label=above left:$c_2$] at (c2) {};
				\node[vtx,label=above left:$c_3$] at (c3) {};
				\node[vtx,label=right:$v$] at (v) {};
			\end{tikzpicture}
			\caption{Extension across $S=P\cap B$ in Theorem~\ref{thm:chordal}. Solid, dashed and dotted lines represent $T_1,T_2,T_3$, respectively. For each $i$, the line from $r$ to $c_i$ represents the path $rT_ic_i$. The three edges $c_iv$ extend these paths to $v$, and the resulting paths meet pairwise only at $r$ and $v$. Other clique vertices and edges are omitted.}
			\label{fig:cliqueextension}
		\end{figure}
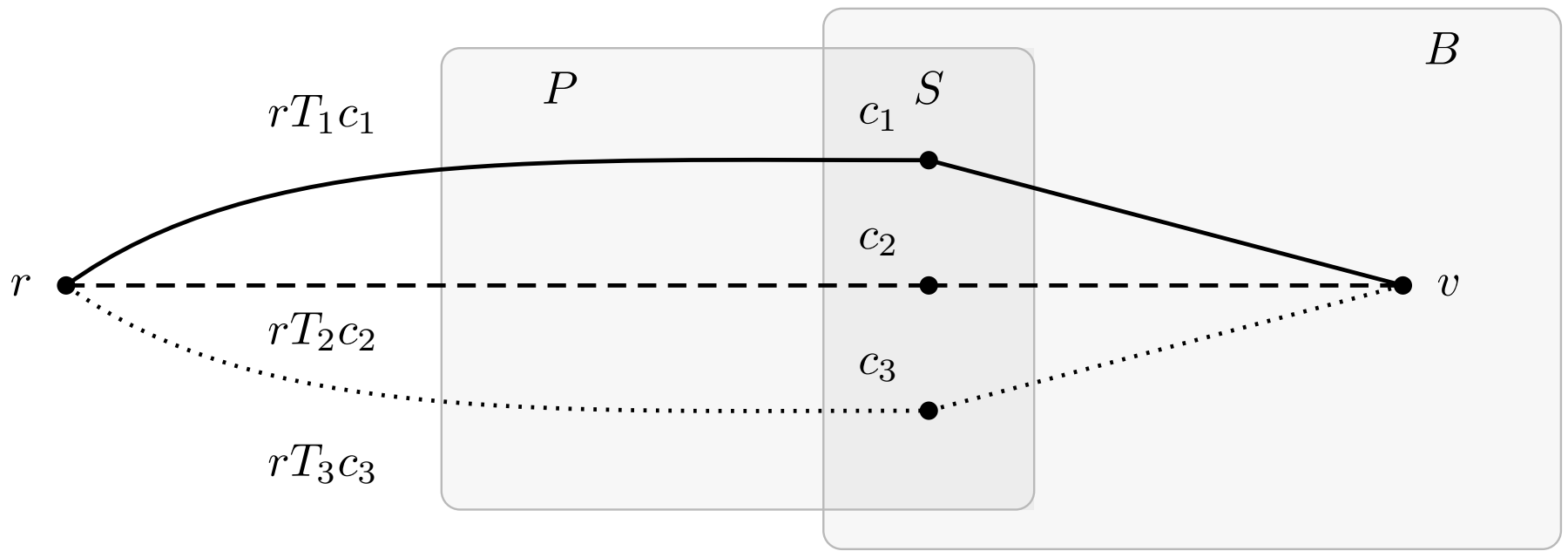
		
		For each $i$, take the union of the trees constructed at all stages, retaining the notation $T_i$. Every vertex of $G$ appears at a finite stage and retains its path to $r$. Hence $T_i$ is spanning and connected. The edges of any cycle in $T_i$ would all occur by some stage, contrary to acyclicity at that stage. The root paths to each fixed vertex are those at its first stage, so the three trees remain independent. The orders $\prec_i$ also extend to the unions. Finally, each $T_i$ is locally finite because it is a subgraph of $G$.
		
		The following stabilization statement is the central lemma in the proof of end-faithfulness. It shows how the rule of choosing the least available attachment vertex forces all vertices of a sufficiently distant bag on a clique-tree ray to have the same ancestor at any prescribed depth.
		
		\begin{claim}\label{cl:descendants}
			Fix $i\in\{1,2,3\}$, an integer $d\geq1$, and a ray $B_0B_1\ldots$ of $\mathcal D$ starting at its root. There are a vertex $z$ of depth $d$ in $T_i$ and an integer $N$ such that every vertex of $B_n$ is a descendant of $z$ in $T_i$ whenever $n\geq N$.
		\end{claim}
		
		\begin{proof}
			The vertices of $T_i$ at each fixed depth form a finite set, by induction on the depth and local finiteness. Thus there are only finitely many vertices at depth less than $d$. Each belongs to finitely many bags, so all sufficiently distant bags on the given ray avoid them. In these bags, assign to each vertex its unique ancestor at depth $d$ as a label. There are finitely many possible labels. Order them by $\prec_i$. The definition of this order implies that every vertex with an earlier label precedes every vertex with a later label.
			
			For all sufficiently large $n$, a vertex of $B_{n+1}\cap B_n$ retains its label, while each vertex of $B_{n+1}\setminus B_n$ was introduced when $B_{n+1}$ was processed, by Equation~\ref{eq:cliqueattachment}. In $T_i$, every such new vertex is attached specifically to $b_i(B_{n+1})\in B_{n+1}\cap B_n$ and therefore receives the label of $b_i(B_{n+1})$. Thus a label absent from $B_n$ cannot occur in $B_{n+1}$. The nonempty sets of labels occurring in the bags $B_n$ consequently form a decreasing sequence. Since there are finitely many possible labels, these sets are eventually constant. Restrict attention to such a tail and let $\lambda$ be the least label in the constant set.
			
			Vertices with label $\lambda$ are introduced arbitrarily far along this tail. Otherwise every occurrence of $\lambda$ in its bags would be accounted for by finitely many vertices in a finite initial union of bags. Each of those vertices belongs to only finitely many bags, so eventually the label would be absent, a contradiction. Choose a sufficiently large $m$ at which a new vertex with label $\lambda$ is introduced. Its parent is $b_i(B_m)$, which consequently also has label $\lambda$.
			
			We show that $b_i(B_n)$ has label $\lambda$ for every $n\geq m$. Suppose otherwise, and let $n>m$ be the first index for which its label differs. The vertex $b_i(B_{n-1})$ was not retained, since retention would preserve its label. When $c_i=b_i(B_n)$ was selected from $S=B_{n-1}\cap B_n$, every vertex of label $\lambda$ in $S$ had already been selected for another index. Indeed, an unselected vertex with that label would precede every available vertex of any other label. Since all new vertices of $B_n$ have the label of $b_i(B_n)$, every vertex of $B_n$ with label $\lambda$ is $b_j(B_n)$ for some $j\neq i$.
			
			We now make explicit why $b_i$ cannot subsequently reacquire label $\lambda$. Let $L$ be the set of vertices of $B_n$ with that label, so $|L|\leq2$. Each vertex of $L$ remains the chosen vertex for its index $j\neq i$ for as long as it occurs in successive bags: all surviving chosen vertices are retained before any replacement is selected. Moreover, once a vertex disappears from the ray of bags it cannot reappear, by the connectedness of the bags containing it.
			
			Inductively, suppose that at some stage $k\geq n$ the only vertices of $B_k$ with label $\lambda$ are surviving members of $L$, each still chosen for its original index $j\neq i$, and that $b_i(B_k)$ has a different label. Every vertex with label $\lambda$ in $B_k\cap B_{k+1}$ is retained for its other index before the choice for $i$ is made, so none is available to $i$. If $b_i(B_k)$ survives, it is retained and still has a different label; otherwise its replacement also has a different label. Hence $b_i(B_{k+1})$ cannot reacquire $\lambda$. Every new vertex of $B_{k+1}\setminus B_k$ is attached in $T_i$ to this particular vertex $b_i(B_{k+1})$, so no new vertex receives label $\lambda$ either. This proves the induction step.
			
			Thus the only vertices that can carry label $\lambda$ in subsequent bags are the at most two members of $L$. Each belongs to finitely many bags, so eventually all of them disappear. Once they have disappeared there is no remaining source from which the label can be recreated, since a new vertex always inherits the label of its attachment vertex in the preceding bag. This contradicts the constancy of the label set.
			
			Figure~\ref{fig:labelstabilization} illustrates this contradiction.
			
			\begin{figure}[htbp]
				\centering
				\begin{tikzpicture}[
					x=1cm,y=1cm,
					font=\fontsize{10}{12}\selectfont,
					vtx/.style={circle,fill=black,inner sep=1.4pt},
					other label/.style={circle,draw=black,fill=white,inner sep=1.8pt,line width=.9pt},
					selected/.style={line width=.9pt},
					omitted/.style={draw=gray!60,dashed,line width=.45pt},
					oldregion/.style={draw=gray!55,fill=gray!6,line width=.45pt,rounded corners=4pt},
					bag transition/.style={-{Stealth[length=3pt,width=4pt]},line width=.6pt},
					every node/.append style={outer sep=1pt}
					]
					\foreach \x in {0,3.3,6.6,9.9}{
						\path[oldregion] (\x,0) rectangle ++(2.7,3.4);
					}
					\node at (1.35,3.75) {$B_{n-1}$};
					\node at (4.65,3.75) {$B_n$};
					\node at (7.95,3.75) {$B_{n+1}$};
					\node at (11.25,3.75) {$B_\ell$};
					\draw[bag transition] (2.77,1.7)--(3.23,1.7);
					\draw[bag transition] (6.07,1.7)--(6.53,1.7);
					\node[inner sep=0pt] at (9.6,1.7) {$\cdots$};
					\node[vtx,label=right:{$x=b_i$}] at (.4,2.9) {};
					\node[vtx,label=right:{$u=b_j$}] at (.4,2.1) {};
					\node[vtx,label=right:{$w=b_k$}] at (.4,1.3) {};
					\node[other label,label=right:$a$] at (.4,.5) {};
					\node[other label,label=right:{$a=b_i$}] at (3.7,2.9) {};
					\node[vtx,label=right:{$u=b_j$}] at (3.7,2.1) {};
					\node[vtx,label=right:{$w=b_k$}] at (3.7,1.3) {};
					\node[other label,label=right:{$y$ (new)}] at (3.7,.5) {};
					\node[other label,label=right:{$a=b_i$}] at (7,2.9) {};
					\node[other label,label=right:{$y=b_j$}] at (7,2.1) {};
					\node[vtx,label=right:{$w=b_k$}] at (7,1.3) {};
					\node[other label,label=right:{new vertex}] at (7,.5) {};
					\node[other label,label=right:$b_i$] at (10.3,2.9) {};
					\node[other label,label=right:$b_j$] at (10.3,2.1) {};
					\node[other label,label=right:$b_k$] at (10.3,1.3) {};
					\node[other label,label=right:{$\cdots$}] at (10.3,.5) {};
					\node[vtx,label=right:{label $\lambda$}] at (2.05,-.65) {};
					\node[other label,label=right:{label different from $\lambda$}]
					at (6.1,-.65) {};
				\end{tikzpicture}
				\caption{The contradiction in Claim~\ref{cl:descendants} when $b_i(B_n)$ does not have label $\lambda$. In the example shown, $y$ first appears in $B_n$ and is chosen as $b_j$ in $B_{n+1}$. Surviving vertices labelled $\lambda$ remain chosen for other indices, so no new vertex receives $\lambda$. For sufficiently large $\ell$, all vertices labelled $\lambda$ have disappeared from $B_\ell$, contradicting the fixed label set. Each $b_j$ refers to its panel's bag.}
				\label{fig:labelstabilization}
			\end{figure}
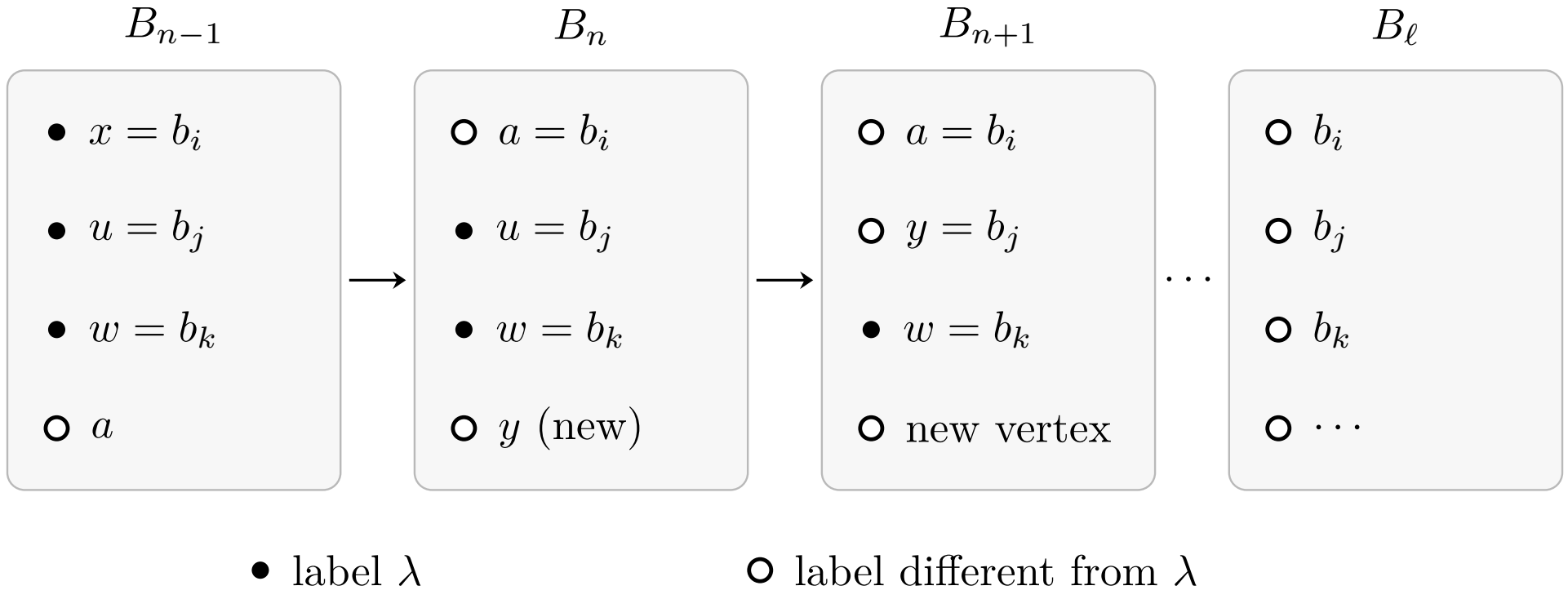
			
			We have proved that $b_i(B_n)$ has label $\lambda$ for all $n\geq m$. Every vertex introduced in a bag $B_n$ with $n>m$ has that label as well. Vertices with other labels can therefore occur only among the finitely many vertices already in $B_m$, and each of these eventually disappears. All sufficiently distant bags consequently have the single label $\lambda$. Taking $z=\lambda$ proves Claim~\ref{cl:descendants}.
		\end{proof}
		
		We describe how a ray of $G$ determines a ray of $\mathcal D$. For each integer $h\geq0$, let $\mathcal D^{\leq h}$ be the finite subtree induced by the bags at distance at most $h$ from $B_0$, and let $W_h$ be the union of those bags. This union is finite. If $v\notin W_h$, the bags containing $v$ form a connected subtree disjoint from $\mathcal D^{\leq h}$ and hence lie in a single component of $\mathcal D-\mathcal D^{\leq h}$. If two vertices outside $W_h$ are adjacent in $G$, they lie together in a bag, so their associated components are equal. It follows that each component of $G-W_h$ is associated with a single component of $\mathcal D-\mathcal D^{\leq h}$.
		
		Let $R$ be a ray of $G$. Its tail outside $W_h$ is connected and thus determines one such component $Q_h$ of $\mathcal D-\mathcal D^{\leq h}$. These components are nested: $Q_{h+1}\subseteq Q_h$, as follows by considering a common tail outside $W_h\cup W_{h+1}$. Every component of $\mathcal D-\mathcal D^{\leq h}$ has a unique vertex at distance $h+1$ from $B_0$. Taking this vertex from each $Q_h$ gives a ray $B_0B_1\ldots$ of $\mathcal D$. If two rays of $G$ belong to the same end, they have tails in the same component of $G-W_h$ for every $h$, so they give the same $Q_h$ and the same ray of $\mathcal D$.
		
		Suppose now that two distinct rays $R,R'$ of $T_i$ starting at $r$ belong to the same end of $G$. They determine the same ray $B_0B_1\ldots$ of $\mathcal D$. Since distinct rays in a tree starting at the same vertex cannot meet again after separating, choose $d\geq1$ larger than the depth of their last common vertex. The vertices of $R$ and $R'$ at depth at least $d$ have different ancestors at depth $d$.
		
		For each $n\geq1$, let $S_n=B_{n-1}\cap B_n$, and let $X_n$ be the union of the bags in the component of $\mathcal D-B_{n-1}B_n$ containing $B_n$. The tails of both $R$ and $R'$ lie in $X_n\setminus S_n$, by their construction of the same component $Q_{n-1}$. For sufficiently large $n$, we also have $r\notin X_n$, since $r$ belongs to only finitely many bags. Proposition~\ref{prop:decompositionseparation} therefore forces both rays, which start at $r$, to meet $S_n$. By Claim~\ref{cl:descendants}, for all sufficiently large $n$ every vertex of $B_n$, and hence of $S_n$, has the same ancestor at depth $d$. In particular, their depths are at least $d$. This contradicts the fact that the vertices of $R$ and $R'$ at those depths have different ancestors at depth $d$.
		
		Thus no two rays of $T_i$ starting at $r$ represent the same end of $G$. Proposition~\ref{prop:endrepresentation} supplies a ray of $T_i$ starting at $r$ in every end of $G$, since $G$ is connected and locally finite and $T_i$ is a spanning tree. Hence $T_i$ is end-faithful in the sense of Definition~\ref{def:faithful}, for each $i$.
	\end{proof}
	
	The local finiteness assumption in Theorem~\ref{thm:chordal} cannot be omitted, even when only one vertex has infinite degree.
	
	\begin{example}\label{ex:squarefan}
		Let $G$ have vertex set $\{z\}\cup\{v_n:n\geq0\}$ and edge set
		\[
		E(G)=\{zv_n:n\geq0\}\cup
		\{v_iv_j:i,j\geq0,\ 1\leq |i-j|\leq2\}.
		\]
		Figure~\ref{fig:squarefan} shows the graph and the two-vertex separator used below after deleting $z$. Notice that $G$ is chordal, $3$-connected and has a unique end, but has no three independent end-faithful spanning trees rooted at any vertex. Only $z$ has infinite degree.

		\begin{figure}[htbp]
			\centering
			\begin{tikzpicture}[
				x=1cm,y=1cm,
				font=\fontsize{10}{12}\selectfont,
				vtx/.style={circle,fill=black,inner sep=1.4pt},
				selected/.style={draw=black,line width=.9pt},
				omitted/.style={draw=gray!60,dashed,line width=.45pt},
				oldregion/.style={fill=gray!14,rounded corners=4pt},
				every node/.append style={outer sep=1pt}
				]
				\path[oldregion] (5.25,.4) rectangle (7.85,1.65);
				\coordinate (z) at (5.7,3.8);
				\foreach \i in {0,...,6}{
					\coordinate (v\i) at ({.6+1.7*\i},1.2);
					\draw[omitted] (z)--(v\i);
				}
				\foreach \i in {0,...,5}{
					\pgfmathtruncatemacro{\j}{\i+1}
					\draw[selected] (v\i)--(v\j);
				}
				\foreach \i in {0,...,4}{
					\pgfmathtruncatemacro{\j}{\i+2}
					\draw[selected] (v\i) .. controls +(.7,-1.45) and +(-.7,-1.45) .. (v\j);
				}
				\begin{scope}
					\clip (0,-.3) rectangle (11.35,4);
					\draw[omitted] (z)--(12.5,1.2);
					\draw[selected] (v6)--(12.5,1.2);
					\draw[selected] (v5) .. controls +(.7,-1.45) and +(-.7,-1.45) .. (12.5,1.2);
					\draw[selected] (v6) .. controls +(.7,-1.45) and +(-.7,-1.45) .. (14.2,1.2);
				\end{scope}
				\node at (11.75,1.2) {$\cdots$};
				\foreach \i in {0,...,6}{
					\node[vtx] at (v\i) {};
					\node at ({.6+1.7*\i},-.35) {$v_{\i}$};
				}
				\node[vtx,label=above:$z$] at (z) {};
			\end{tikzpicture}
			\caption{The graph of Example~\ref{ex:squarefan}. The black edges $v_iv_{i+1}$ and $v_iv_{i+2}$ and the dashed edges $zv_i$ continue in the same way for all $i\geq0$. After deleting $z$, the shaded pair $\{v_3,v_4\}$ separates the initial segment from the tail. Such pairs exclude three vertex-disjoint rays.}
			\label{fig:squarefan}
		\end{figure}
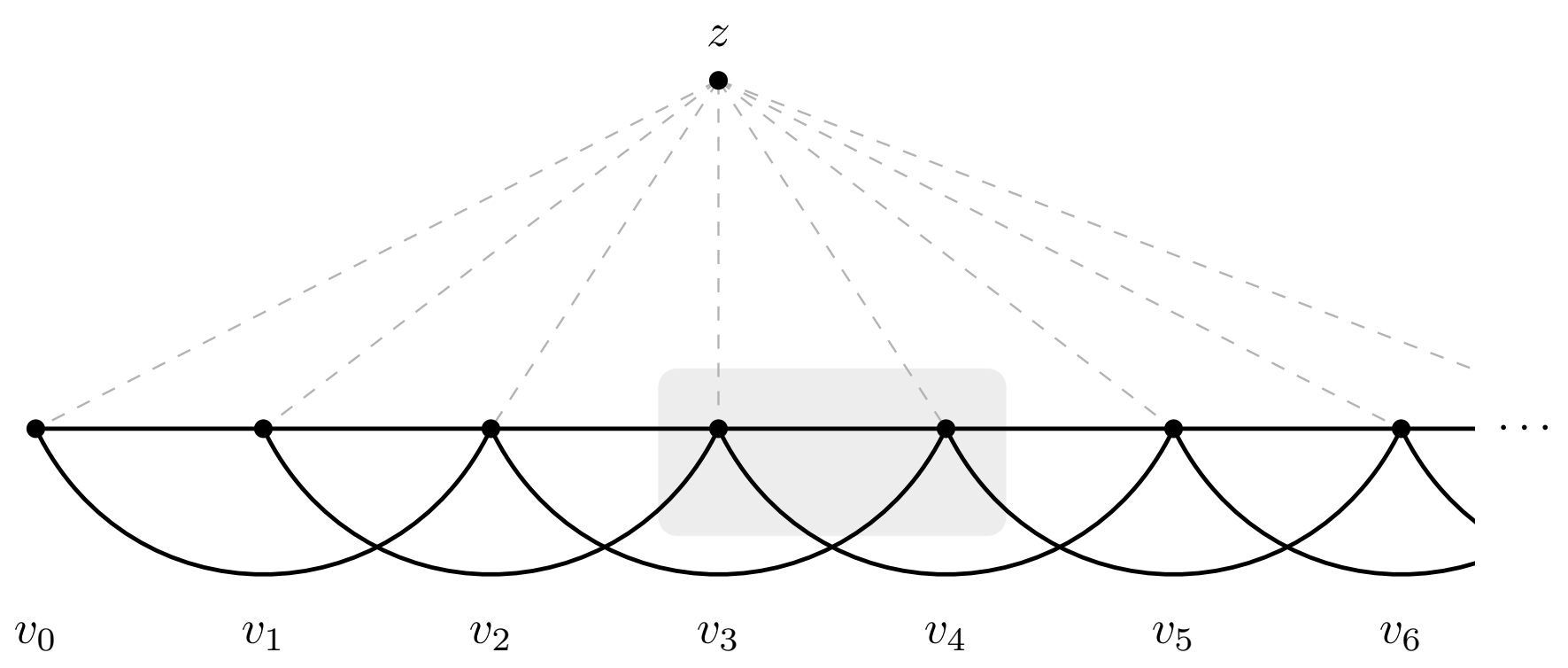
		
		To verify the connectivity assertion, delete at most two vertices. If $z$ remains, it joins all other remaining vertices. If $z$ is deleted, at most one $v_i$ is also deleted. If $v_0$ is deleted, the ray $v_1v_2\ldots$ spans the remaining graph, while a deleted $v_i$ with $i\geq1$ is bypassed by the edge $v_{i-1}v_{i+1}$. Thus the remaining graph is connected, and $G$ is $3$-connected. Deleting $z,v_1,v_2$ isolates $v_0$ and leaves other vertices, so its connectivity is exactly three.
		
		Every cycle of length at least four containing $z$ has a chord incident with $z$. In a cycle of length at least four avoiding $z$, let $v_i$ have the smallest index among its vertices. Its two neighbours on the cycle must be $v_{i+1}$ and $v_{i+2}$, since those are the only neighbours of $v_i$ with larger index. They are adjacent and are not consecutive on this cycle, so their edge is a chord. Hence $G$ is chordal.
		
		Every ray of $G$ has a tail avoiding $z$, and its vertex indices on that tail are unbounded. Given a finite vertex set $X$, choose $n$ larger than every index occurring in $X$. A ray visits each of the finitely many vertices $z,v_0,\ldots,v_{n-1}$ at most once, so it has a tail in $\{v_j:j\geq n\}$. These vertices induce a connected subgraph disjoint from $X$. All rays therefore belong to the same end, by Definition~\ref{def:ends}. At least one ray exists, namely $v_0v_1\ldots$, so this end is unique.
		
		There are no three pairwise vertex-disjoint rays in $G$. If there were, take tails avoiding $z$ and choose $n$ greater than the indices of their three initial vertices. Every tail has unbounded indices and must therefore meet $\{v_n,v_{n+1}\}$: an edge of $G-z$ changes the index by at most two, so a path from an index less than $n$ to one greater than $n+1$ cannot avoid this set. Three disjoint tails cannot all meet a set of two vertices.
		
		If three independent end-faithful spanning trees existed, each would contain a ray starting at the common root in the unique end of $G$. The observation in Example~\ref{ex:fan} would give pairwise vertex-disjoint tails of these three rays, a contradiction. This argument applies to every choice of root.
	\end{example}

	Theorems~\ref{thm:three} and~\ref{thm:chordal} leave the following questions.
	
	\begin{question}\label{ques:three}
		Does every infinite $3$-connected graph have three independent spanning trees rooted at any prescribed vertex?
	\end{question}
	
	\begin{question}\label{ques:threefaithful}
		Does every locally finite $3$-connected graph have three independent end-faithful spanning trees rooted at any prescribed vertex?
	\end{question}
	
	\section*{Acknowledgements}
	This study was financed in part by the Coordena\c{c}\~ao de Aperfei\c{c}oamento de Pessoal de N\'ivel Superior - Brasil (CAPES) - Finance Code 001.
	
	\section*{Declaration of AI-assisted use}
	
The author wrote the manuscript, provided mathematical input during the
development of the proofs of Theorems~\ref{thm:two}, \ref{thm:three},
\ref{thm:faithful} and~\ref{thm:chordal}, and constructed the
counterexamples in Examples~\ref{ex:fan} and~\ref{ex:squarefan}.
ChatGPT contributed substantially to the development of the proofs
and was used to discuss the exposition. Claude was also used to review
the proofs for possible gaps, check the consistency of the notation,
and suggest improvements to the writing.

The author rewrote all AI-assisted arguments, verified every step
of the proofs, and takes full responsibility for the final manuscript.

	\Addresses
	

\begin{thebibliography}{99}
		
		\bibitem{BK}
		C.~B\"urger and J.~Kurkofka.
		End-faithful spanning trees in graphs without normal spanning trees.
		\emph{J. Graph Theory}, 101(1):95--105, 2022.
		\href{https://doi.org/10.1002/jgt.22813}
		{doi:\nolinkurl{10.1002/jgt.22813}}.
		
		\bibitem{CWF}
		B.~Cheng, D.~Wang, and J.~Fan.
		Independent spanning trees in networks: A survey.
		\emph{ACM Comput. Surv.}, 55(14s):Article 335, 2023.
		\href{https://doi.org/10.1145/3591110}
		{doi:\nolinkurl{10.1145/3591110}}.
		
		\bibitem{CM}
		J.~Cheriyan and S.~N.~Maheshwari.
		Finding nonseparating induced cycles and independent spanning trees
		in $3$-connected graphs.
		\emph{J. Algorithms}, 9(4):507--537, 1988.
		\href{https://doi.org/10.1016/0196-6774(88)90015-6}
		{doi:\nolinkurl{10.1016/0196-6774(88)90015-6}}.
		
		\bibitem{CLY}
		S.~Curran, O.~Lee, and X.~Yu.
		Finding four independent trees.
		\emph{SIAM J. Comput.}, 35(5):1023--1058, 2006.
		\href{https://doi.org/10.1137/S0097539703436734}
		{doi:\nolinkurl{10.1137/S0097539703436734}}.
		
		\bibitem{Diestel}
		R.~Diestel.
		\emph{Graph Theory}, volume 173 of Graduate Texts in Mathematics.
		Springer, Berlin--Heidelberg, 6th edition, 2025.
		\href{https://doi.org/10.1007/978-3-662-70107-2}
		{doi:\nolinkurl{10.1007/978-3-662-70107-2}}.
		
		\bibitem{DiestelEnds}
		R.~Diestel.
		On end-faithful spanning trees in infinite graphs.
		\emph{Math. Proc. Cambridge Philos. Soc.}, 107(3):461--473, 1990.
		\href{https://doi.org/10.1017/S0305004100068742}
		{doi:\nolinkurl{10.1017/S0305004100068742}}.
		
		\bibitem{DK}
		R.~Diestel and D.~K\"uhn.
		On infinite cycles I.
		\emph{Combinatorica}, 24(1):69--89, 2004.
		\href{https://doi.org/10.1007/s00493-004-0005-z}
		{doi:\nolinkurl{10.1007/s00493-004-0005-z}}.
		
		\bibitem{DFKPY}
		N.~Dragani\'c, K.~Frankston, M.~Krivelevich, A.~Pokrovskiy,
		and L.~Yepremyan.
		On independent spanning trees in random graphs.
		In \emph{Proceedings of the 2026 Annual ACM-SIAM Symposium
			on Discrete Algorithms (SODA)}, pages 4096--4104.
		SIAM, 2026.
		\href{https://doi.org/10.1137/1.9781611978971.151}
		{doi:\nolinkurl{10.1137/1.9781611978971.151}}.
		
		\bibitem{Freudenthal}
		H.~Freudenthal.
		\"Uber die Enden topologischer R\"aume und Gruppen.
		\emph{Math. Z.}, 33:692--713, 1931.
		\href{https://doi.org/10.1007/BF01174375}
		{doi:\nolinkurl{10.1007/BF01174375}}.
		
		\bibitem{Gavril72}
		F.~Gavril.
		Algorithms for minimum coloring, maximum clique, minimum covering
		by cliques, and maximum independent set of a chordal graph.
		\emph{SIAM J. Comput.}, 1(2):180--187, 1972.
		\href{https://doi.org/10.1137/0201013}
		{doi:\nolinkurl{10.1137/0201013}}.
		
		\bibitem{Gavril74}
		F.~Gavril.
		The intersection graphs of subtrees in trees are exactly the chordal graphs.
		\emph{J. Combin. Theory Ser. B}, 16(1):47--56, 1974.
		\href{https://doi.org/10.1016/0095-8956(74)90094-X}
		{doi:\nolinkurl{10.1016/0095-8956(74)90094-X}}.
		
		\bibitem{Halin}
		R.~Halin.
		\"Uber unendliche Wege in Graphen.
		\emph{Math. Ann.}, 157:125--137, 1964.
		\href{https://doi.org/10.1007/BF01362670}
		{doi:\nolinkurl{10.1007/BF01362670}}.
		
		\bibitem{HTL}
		C.~Hofer-Temmel and F.~Lehner.
		Clique trees of infinite locally finite chordal graphs.
		\emph{Electron. J. Combin.}, 25(2):\#P2.9, 2018.
		\href{https://doi.org/10.37236/3928}
		{doi:\nolinkurl{10.37236/3928}}.
		
		\bibitem{HLMPW}
		L.~Hollom, L.~Lichev, A.~Mond, J.~Portier, and Y.~Wang.
		Approximate Itai--Zehavi conjecture for random graphs.
		\emph{Random Structures Algorithms}, 69(2):e70095, 2026.
		\href{https://doi.org/10.1002/rsa.70095}
		{doi:\nolinkurl{10.1002/rsa.70095}}.
		
		\bibitem{IR}
		A.~Itai and M.~Rodeh.
		The multi-tree approach to reliability in distributed networks.
		\emph{Inform. and Comput.}, 79(1):43--59, 1988.
		\href{https://doi.org/10.1016/0890-5401(88)90016-8}
		{doi:\nolinkurl{10.1016/0890-5401(88)90016-8}}.
		
		\bibitem{ST}
		P.~D.~Seymour and R.~Thomas.
		An end-faithful spanning tree counterexample.
		\emph{Proc. Amer. Math. Soc.}, 113(4):1163--1171, 1991.
		\href{https://doi.org/10.1090/S0002-9939-1991-1045600-8}
		{doi:\nolinkurl{10.1090/S0002-9939-1991-1045600-8}}.
		
		\bibitem{Siran}
		J.~\v{S}ir\'a\v{n}.
		End-faithful forests and spanning trees in infinite graphs.
		\emph{Discrete Math.}, 95(1--3):331--340, 1991.
		\href{https://doi.org/10.1016/0012-365X(91)90345-3}
		{doi:\nolinkurl{10.1016/0012-365X(91)90345-3}}.
		
		\bibitem{ZI}
		A.~Zehavi and A.~Itai.
		Three tree-paths.
		\emph{J. Graph Theory}, 13(2):175--188, 1989.
		\href{https://doi.org/10.1002/jgt.3190130205}
		{doi:\nolinkurl{10.1002/jgt.3190130205}}.
		
	\end{thebibliography}
\end{document}